\documentstyle[12pt]{article}
\newtheorem{thm}{Theorem}
\newtheorem{prop}[thm]{Proposition}

\newtheorem{thm-defi}[thm]{Theorem/Definition}
\newtheorem{example}[thm]{Example}
\newtheorem{cor}[thm]{Corollary}

\newtheorem{new-lemma}[thm]{Lemma}
\newtheorem{defi}[thm]{Definition}
\newtheorem{rem}[thm]{Remark}

\newtheorem{condition}[thm]{Condition}

\newcommand{\A}{{\cal A}}

\newcommand{\E}{{\cal E}}
\newcommand{\F}{{\cal F}}

\renewcommand{\H}{{\cal H}}

\newcommand{\M}{{\cal M}}

\newcommand{\U}{{\cal U}}

\renewcommand{\P}{{\cal P}}
\newcommand{\PP}{{\Bbb P}}

\newcommand{\Integers}{{\Bbb Z}}
\newcommand{\ComplexNumbers}{{\Bbb C}}
\newcommand{\RationalNumbers}{{\Bbb Q}}

\newcommand{\IsomRightArrow}{\stackrel{\cong}{\rightarrow}}
\newcommand{\LongIsomRightArrow}{\stackrel{\cong}{\longrightarrow}}

\newcommand{\LongRightArrowOf}[1]{\stackrel{#1}{\longrightarrow}}

\newcommand{\StructureSheaf}[1]{{\cal O}_{#1}}
\newcommand{\EndProof}{\hfill  $\Box$}
\newcommand{\restricted}[2]{#1_{\mid_{#2}}}

\newcommand{\rank}{\rm rank}

\newcommand{\Ext}{{\rm Ext}}

\newcommand{\Aut}{{\rm Aut}}

\newcommand{\SheafTor}{{\cal T}or}

\newcommand{\Choose}[2]
{\left(\!\!\begin{array}{c}#1\\#2\end{array}\!\!\right)}

\input amssymb.sty

\begin{document}
%\rightline{\today}
\begin{center}
\begin{Large}
{\bf 
\noindent
Integral generators for the cohomology ring
of moduli spaces of sheaves over Poisson surfaces
%Generators for the cohomology ring, with integral coefficients, 
%of Hilbert schemes of Poisson surfaces
}
\end{Large}
\\
Eyal Markman
%University of Massachusetts, Amherst, MA 01003, 
%E-mail: markman@math.umass.edu
\end{center}

\begin{abstract}
Let $\M$ be a smooth and compact moduli space of stable coherent sheaves
on a projective surface $S$ with an effective (or trivial) anti-canonical 
line bundle. We find generators for the cohomology ring of 
$\M$, with integral coefficients. When $S$ is simply connected
and a universal sheaf $\E$ exists over $S\times \M$, then 
its class $[\E]$ admits a K\"{u}nneth decomposition as a class 
in the tensor product
$K_{top}^0(S)\otimes K_{top}^0(\M)$ of the topological $K$-rings.  The 
generators are the Chern classes of the K\"{u}nneth factors of $[\E]$ in
$K_{top}^0(\M)$. The general case is similar.
\end{abstract}

{\scriptsize 
\tableofcontents
} 

%***************************************************************************
% Introduction 
%***************************************************************************
\section{Introduction}
\label{sec-introduction}

Let $S$ be a smooth connected projective symplectic or Poisson surface. 
If symplectic, $S$ is either a $K3$ or an abelian surface. 
Non-symplectic projective Poisson surfaces include 
all minimal rational surfaces, all del-Pezzo surfaces, and 
certain ruled surfaces over projective curves of arbitrary genus. 
Poisson surfaces are classified in \cite{bartocci-macri}.
%Note, that the cohomology groups $H^i(S,\Integers)$ are torsion free,
%for any projective symplectic or Poisson surface. 

Given an algebraic variety $X$, 
denote by $K_{alg}^0(X)$ the Grothendieck $K$-ring of algebraic vector 
bundles on $X$ and by $K_{top}^0(X)$ its topological analogue.
Let $v$ be a class in $K^0_{top}(S)$ of rank
$r\geq 0$ with $c_1(v)$ of Hodge-type $(1,1)$.
Assume that $v$ is primitive, i.e., $v$ is not a multiple of
another class in $K_{top}^0(S)$ by an integer $\geq 1$. 
If $r=0$, assume for simplicity, that the anti-canonical 
line-bundle is either trivial or ample (the assumption is relaxed in
Condition \ref{cond-proper-subsheaf}). 
Given an ample line bundle  $H$ on $S$, 
denote by $\M:=\M_H(v)$ the moduli space of $H$-stable sheaves on $S$ 
with class $v$. We use stability in the sense of Gieseker, Maruyama,  
and Simpson (Definition \ref{def-stability}). 
For a generic choice of an ample line bundle  $H$ on $S$,
called $v$-{\em generic\/} in Definition \ref{def-v-suitable}, 
the moduli space $\M_H(v)$ is either empty, or smooth of the
expected dimension, connected, and 
projective, and it admits a holomorphic symplectic or Poisson structure
(section \ref{sec-stable-sheaves-and-modui-spaces}). 
The expected dimension of $\M_H(v)$ is $\epsilon-\chi(v^\vee\cup v)$, 
where $v^\vee$ is the class dual to $v$, $\cup$ is the product operation
in $K_{top}^0(S)$,
$\epsilon=2$ if $S$ is symplectic, $\epsilon=1$ if
$S$ is non-symplectic but Poisson, and 
$\chi$ is the Euler characteristic defined in section \ref{sec-notation}.

A {\em universal sheaf\/}  is a coherent sheaf $\E$ over $S\times \M$, 
flat over $\M$, 
whose restriction to $S\times \{m\}$,  $m\in\M$, is isomorphic 
to the sheaf $E_m$ on $S$ in the isomorphism class $m$.
The universal sheaf is canonical only up to tensorization by
the pull-back of a line-bundle on $\M$. 
A universal sheaf exists, if there exists a class $x$ in
$K_{alg}^0(S)$ satisfying $\chi(x\cup v)=1$
(\cite{mukai-hodge} or section \ref{sec-universal-class} below). 
Otherwise, there is a weaker notion of a 
{\em twisted} universal sheaf, denoted also by $\E$, 
where the twisting is encoded by a class $\theta$ in \v{C}ech cohomology
$H^2(\M,\StructureSheaf{\M}^*)$, in the classical topology
(Definition \ref{def-twisted-sheaves}). 
For a triple $S$, $v$, $H$, as above, 
the class $\theta$ is always topologically trivial;
it maps to $0$ in $H^3(\M,\Integers)$ via
the connecting homomorphism of the exponential sequence
(Lemma \ref{lemma-cohomological-triviality-of-theta}). 
Consequently, $\E$ defines a class $e$ in $K_{top}^0(S\times\M)$,
canonical up to tensorization by 
the pull-back of the class of a topological line-bundle on $\M$
(Definition \ref{def-e-v}). 

When $H^1(S,\Integers)$ does not vanish, 
we will need the odd $K$-groups of $S$ and $\M$ as well. 
Given a complex algebraic variety $X$, let $K_{top}^1(X)$ be its odd
$K$-group and 
$K_{top}^*(X):=K_{top}^0(X)\oplus K_{top}^1(X)$ its $K$-ring
%and the Chern character ring isomorphism
%$ch:K_{top}^*(X)\otimes \RationalNumbers\rightarrow H^*(X,\RationalNumbers)$
(\cite{atiyah-book,karoubi} and section \ref{sec-non-simply-connected}).
The $K$-ring $K^*(S)$ of a projective Poisson surface 
is known to be torsion free. As a consequence, 
the K\"{u}nneth Theorem yields an isomorphism
\begin{equation}
\label{eq-exterior-product-from-cartesian-square-of-moduli-space}
[K^0_{top}(S)\otimes K^0_{top}(\M)] \ \ \oplus \ \ 
[K^1_{top}(S)\otimes K^1_{top}(\M)] \ \ \ \longrightarrow \ \ \ 
K^0_{top}(S\times \M)
\end{equation}
(\cite{atiyah-book}, Corollary 2.7.15). 
Choose a basis $\{x_1, \dots, x_n\}$ of $K^*_{top}(S)$,
which is a union of bases of
the summands $K^0_{top}(S)$ and $K^1_{top}(S)$. 
We get the K\"{u}nneth decomposition
\begin{equation}
\label{introduction-eq-e-i}
e \ \ = \ \ 
\sum_{i=1}^n x_i\otimes e_i,
\end{equation}
of the class 
$e \in K^0_{top}(S\times \M_H(v))$ of the universal sheaf, 
where each class $e_i$ is either in $K^0_{top}(\M_H(v))$ or in 
$K^1_{top}(\M_H(v))$. The Chern classes 
$c_i(y)\in H^{2i}(\M,\Integers)$, for odd classes 
$y\in K^1_{top}(\M)$ and $i\geq 1/2$ a half-integer, 
are introduced in Definition \ref{def-odd-chern-classes}.

The main result of this paper is the following.

\begin{thm}
\label{thm-introduction-integral-generators}
Let $S$ be a projective Poisson surface, $v\in K_{top}^0(S)$ as above, 
and $H$ a $v$-generic polarization.
\begin{enumerate}
\item
\label{thm-item-chern-classes-of-kunneth-factors-generate}
The cohomology ring $H^*(\M_H(v),\Integers)$ is generated by 
the Chern classes $c_j(e_i)$, of the 
K\"{u}nneth factors $e_i\in K^*_{top}\M_H(v)$,
which are given in equation
(\ref{introduction-eq-e-i}). 
\item
\label{thm-item-torsion-free}
The cohomology groups $H^i(\M_H(v),\Integers)$ are torsion free for all $i$.
\item
\label{thm-item-no-odd-cohomology}
If $H^1(S,\Integers)=0$, then $H^i(\M_H(v),\Integers)$
vanishes for odd $i$.
\end{enumerate}
\end{thm}

The Theorem is a summary of Propositions 
\ref{prop-integral-generators} and
\ref{prop-integral-generators-non-simply-connected-surface} and 
Corollary \ref{cor-generators-in-absence-of-universal-sheaf}.
Part
% \ref{thm-item-torsion-free} and 
\ref{thm-item-no-odd-cohomology} 
of the Theorem is an immediate consequence of part
\ref{thm-item-chern-classes-of-kunneth-factors-generate}. 

In a sequel to this paper we 
apply Theorem \ref{thm-introduction-integral-generators} to
the study of the Hilbert schemes $S^{[n]}$, $n\geq 2$, of length $n$ 
zero-dimensional subschemes of a $K3$ surface $S$
\cite{markman-sequel}. 
The Hilbert scheme $S^{[n]}$ has complex deformations, 
which are not Hilbert schemes on any $K3$ surface. 
%$H^2(S^{[n]},\Integers)$ has a natural symmetric non-degenerate 
%bilinear pairing. 
We determine the subgroup ${\rm Mon}^2$ of 
$\Aut[H^2(S^{[n]},\Integers)]$, generated by monodromy operators
of families of deformations of $S^{[n]}$, consisting of smooth 
hyperk\"{a}hler varieties.
We find an arithmetic obstruction for elements of 
$\Aut[H^2(S^{[n]},\Integers)]$ to extend to monodromy operators of
the full integral cohomology ring $H^*(S^{[n]},\Integers)$.
${\rm Mon}^2$ turns out to be
smaller than expected, when $n-1$ is not a prime power. 
As a consequence, we get that the weight $2$ Hodge structure, of
hyperk\"{a}hler deformations $X$ of  $S^{[n]}$,
does not determine the bimeromorphic class of $X$, 
when $n-1$ is not a prime power. 

Denote by $A^*(\M_H(v))$ the Chow ring of $\M_H(v)$. 
The following theorem 
is proven in section 
\ref{sec-a-decomposition-of-the-diagonal-in-the-chow-ring}.

\begin{thm}
\label{thm-introduction-chow-ring-isomorphic-to-cohomology-ring}
Let $S$ be a rational Poisson surface, and $v$, $H$ as 
in Theorem \ref{thm-introduction-integral-generators}. 
The class of the diagonal in $\M_H(v)\times \M_H(v)$ is in the image of the
exterior-product homomorphism
\[
A^*(\M_H(v))\otimes A^*(\M_H(v)) \ \ \rightarrow \ \ 
A^*(\M_H(v)\times \M_H(v)).
\]
Furthermore, the natural homomorphism
$A^*(\M_H(v))\rightarrow H^*(\M_H(v),\Integers)$ is an isomorphism. 
\end{thm}

%Theorem \ref{thm-introduction-chow-ring-isomorphic-to-cohomology-ring}

{\bf Higgs bundles:}
%We find below generators for the integral cohomology 
%of the moduli space of rank $r$ Higgs bundles of degree $d$ over
%a Riemann surface $\Sigma$ of genus $\geq 2$, assuming
%that $r$ and $d$ are relatively prime (Theorem \ref{thm-higgs}). 
%This 
Let $\Sigma$ be a smooth compact and connected Riemann surface of genus 
$g\geq 2$, and $D$ an effective divisor on $\Sigma$. 
The important special case $D=0$ is included.
A {\em Higgs bundle on $\Sigma$,
with possible poles along $D$}, is a pair $(E,\varphi)$,
consisting of a vector bundle $E$ on $\Sigma$ and a 
$1$-form valued endomorphism $\varphi:E\rightarrow E\otimes K_\Sigma(D)$.
The Higgs bundle is {\em stable}, if any non-zero proper $\varphi$-invariant
subbundle $F\subset E$ satisfies the 
inequality $\deg(F)/\rank(F)< \deg(E)/\rank(E)$. 
The moduli space
$\H_\Sigma(r,d,D)$, of stable rank $r$ Higgs bundles of degree $d$, with
possible poles along $D$,
is a smooth quasi-projective variety \cite{nitsure,simpson}. 
Denote by $f_i$, $i=1,2$, the projections from 
$\Sigma\times \H_\Sigma(r,d,D)$.
Assume that $r$ and $d$ are relatively prime. Then there exists 
over $\Sigma\times \H_\Sigma(r,d,D)$ a universal vector bundle $\E$ and
a universal Higgs field
$\Phi:\E\rightarrow \E\otimes f_1^*K_\Sigma(D)$. 
Choose a basis $\{x_1, x_2, \dots, x_{2g}\}$ of
$K^1_{top}(\Sigma)$, and $\{x_{2g+1}, x_{2g+2}\}$ of 
$K^0_{top}(\Sigma)$.
%let $x_0$ be the class of the trivial line bundle,
%and $x_{2g+1}$ the class of the sky-scraper sheaf of a point. 
Define the K\"{u}nneth factors 
$e_i\in K_{top}^*\left(\H_\Sigma(r,d,D)\right)$,
$1\leq i \leq 2g+2$, of the universal bundle $\E$, as in equation
(\ref{introduction-eq-e-i}). 
%See Definition
%\ref{def-odd-chern-classes} for 
%the Chern classes $c_j(e_i)$, $1\leq i\leq 2g$.

\begin{thm}
\label{thm-higgs}
The cohomology ring $H^*\left(\H_\Sigma(r,d,D),\Integers\right)$ is 
generated by the Chern classes $c_j(e_i)$. 
\end{thm}

The Theorem is proven in section
\ref{sec-higgs}. It
sharpens results of Hausel-Thaddeus 
and the author, where the cohomology was considered with 
rational coefficients \cite{ht,markman-diagonal}.

{\bf Related works:}
Ellingsrud and Str{\o}mme proved Theorems 
\ref{thm-introduction-integral-generators} and 
\ref{thm-introduction-chow-ring-isomorphic-to-cohomology-ring}
when the surface $S$ is the projective plane. 
Beauville found generators, for the cohomology ring with 
{\em rational} coefficients
of moduli spaces as above, when $X$ is a non-symplectic Poisson surface 
%rational or ruled surface
\cite{beauville-diagonal}. 
When $S$ is a K3 or abelian surface, 
generators for the cohomology ring $H^*(\M,\RationalNumbers)$,
with {\em rational} coefficients, 
where found in \cite{markman-diagonal}. 
In case $S$ is an arbitrary projective surface and 
$\M=S^{[n]}$ is the Hilbert scheme, parametrizing
ideals sheaves of length $n$ subscheme, generators for the ring 
$H^*(S^{[n]},\RationalNumbers)$ were found
in \cite{lqw1,lqw2}. 
The ring structure of 
$H^*(S^{[n]},\RationalNumbers)$ was calculated by Lehn and Sorger
for a K3 or abelian surface $S$
\cite{lehn-sorger}. The cohomology ring
was shown to be isomorphic to the orbifold cohomology of the symmetric product
\cite{fantechi-gottsche}. 
The ring structure of $H^*(S^{[n]},\RationalNumbers)$, for any smooth 
projective surface $S$, was determined in
\cite{costello-grojnowski}.

Theorem \ref{thm-introduction-integral-generators}, 
about integral generators for the cohomology of moduli spaces, 
follows easily from the Main Theorems of \cite{beauville-diagonal} and 
\cite{markman-diagonal}, stated as Theorem \ref{thm-diagonal} 
below. Theorem \ref{thm-diagonal} expresses the class of the
diagonal in $\M\times \M$, in terms of the universal sheaf $\E$ over
$S\times \M$. The K\"{u}nneth decomposition 
(\ref{introduction-eq-e-i}) 
of the class of $\E$,
in the topological $K$-group of $S\times \M$, leads to a decomposition
of the class of the diagonal in $\M\times \M$. More precisely, 
the class of the diagonal is Poincare dual to 
the sum of exterior products of Chern classes of the classes 
$e_i\in K_{top}^*(\M)$ given in (\ref{introduction-eq-e-i}). 
Theorem \ref{thm-introduction-integral-generators}
follows easily from the latter decomposition. 
Our decomposition of the diagonal 
is the precise topological analogue of that of \cite{ellingsrud-stromme-p2}, 
replacing $\PP^2$ by a symplectic or Poisson surface $S$, 
and the algebraic $K$-group by the topological one. 

The paper is organized as follows. 
In section \ref{sec-integral-generators} we prove Theorem
\ref{thm-introduction-integral-generators} for 
moduli spaces of sheaves over Poisson surfaces, assuming that 
a universal sheaf exists. 
The case of Poisson surfaces with vanishing odd cohomology
is particularly simple, and is treated first. In section
\ref{sec-twisted-sheaves} we find  integral generators for the cohomology of 
such moduli spaces, when a universal sheaf does not exist. 
In section \ref{sec-higgs} generators are found for the integral cohomology 
ring of the moduli spaces of Higgs bundles.

{\em Acknowledgments}: The work on this paper began, while visiting 
the mathematics department of the University of Lille 1 during June 2003. 
I would like to thank Dimitri Markushevich and Armando Treibich for 
their hospitality. I would like to thank Jonathan Block for
the reference to Lemma
\ref{lemma-reduction-to-torsion-free-case}
of Atiyah. I thank the referees for simplifying the proof of Lemma 
\ref{lemma-invariance-of-c-m-under-tensorization}, for pointing an
error in the original proof of Lemma \ref{lemma-integral-coefficients},
and for additional insightful comments and suggestions.

{\em Note:} Zhenbo Qin and Weiqiang Wang have recently and independently 
obtained bases for the integral cohomology groups (modulo torsion)
of Hilbert schemes of points on a projective surface $X$ 
with vanishing $H^1(X,\StructureSheaf{X})$ and $H^2(X,\StructureSheaf{X})$
\cite{qin-wang}.

%***************************************************************************
% Notation 
%***************************************************************************
\subsection{Notation}
\label{sec-notation}

Given a smooth complex algebraic variety $X$, we denote by 
$K^0_{top}X$ the topological $K$-ring of vector bundles on $X$.
We let $K^0_{alg}X$ be its algebraic analogue and
\[
\alpha \ : \ K^0_{alg}X \ \ \longrightarrow \ \ K^0_{top}X
\]
the natural homomorphism. Given a morphism $f:X\rightarrow Y$,
we denote by $f^!:K^0_{top}Y\rightarrow K^0_{top}X$ the pullback. 
When $f$ is a proper morphism, we denote by
$f_!:K^0_{top}X\rightarrow K^0_{top}Y$ the topological Gysin map
(\cite{bfm}, \cite{karoubi} Proposition IV.5.24). 
We used above the assumption, that $X$ and $Y$ are smooth, 
which enables us to identify $K^0_{top}$ with the $K$-homology groups
$K_0^{top}$, for both spaces. 
Similarly, we identify the Grothendieck group $K_0^{alg}X$, of coherent 
sheaves, with $K^0_{alg}X$, replacing the class of a sheaf by that of
any of its locally free resolutions. 
The algebraic push-forward 
$f_!:K^0_{alg}X\rightarrow K^0_{alg}Y$ 
takes the class of a coherent sheaf $E$ on $X$, 
to the alternating sum $\sum_{i\geq 0} (-1)^i R^i_{f_*}E$ of 
the classes of the higher direct image sheaves on $Y$. 
When $Y$ is a point $\{pt\}$, 
we identify $K^0_{top}(\{pt\})$ with $\Integers$ and $f_!$ is the 
{\em Euler characteristic} 
$\chi:K^0_{top}(X)\rightarrow\Integers$.
The algebraic and topological Gysin maps are compatible via the
equality $f_!\circ\alpha=\alpha\circ f_!$ (see \cite{bfm}).

We denote by $x\cup y$ the product in $K^0_{top}$. 
The Gysin homomorphisms satisfy the projection formula 
\begin{equation}
\label{eq-projection-formula}
f_!(x\cup f^!(y)) \ \ = \ \ f_!(x)\cup y
\end{equation}
both for algebraic and for topological $K$-groups.
%(\cite{atiyah-hirzebruch-rr} and \cite{bfm} page 169). 
The topological Gysin homomorphism satisfies 
the following weak analogue, of the Cohomology and Base Change 
Theorem, applied to a cartesian product instead of a fiber product.
\[
\begin{array}{ccc}
X\times Z & \LongRightArrowOf{f\times \iota} & Y\times Z
\\
\pi_X \ \downarrow \ \hspace{2ex} & & \hspace{2ex} \ \downarrow \ \pi_Y
\\
X & \LongRightArrowOf{f} & Y
\end{array}
\] 
The equivalence 
\begin{equation}
\label{eq-cohomology-and-base-change}
(f\times \iota)_!(\pi_X^!x) \ \ = \ \ 
\pi_Y^!f_!(x)
\end{equation}
holds, 
where $x\in K^0_{top}(X)$, $\iota:Z\rightarrow Z$ is the identity, 
$f\times \iota: X\times Z\rightarrow Y\times Z$ is the product map,
and $\pi_X$, $\pi_Y$ are the projections \cite{atiyah-hirzebruch-rr}. 

$S$ will denote a surface and $\M$ a moduli space of stable sheaves on $S$.
The morphism $f_i$ is the projection from $S\times \M$ on the $i$-th factor,
$i=1,2$. The morphism $p_i$ is the projection from
$\M\times \M$ on the $i$-th factor. 
%The morphism
%$\pi_{ij}$ is the projection from $\M\times S\times \M$ on the
%product of the $i$-th and $j$-th factors. 

%***************************************************************************
% Integral generators via a universal sheaf
%***************************************************************************
\section{Integral generators via a universal sheaf}
\label{sec-integral-generators}
In section \ref{sec-stable-sheaves-and-modui-spaces} the necessary 
background on moduli spaces
of stable sheaves on Poisson surfaces is reviewed.
We recall in section \ref{sec-class-of-diagonal}, 
that the class of the diagonal in $\M\times \M$ can be expressed 
in terms of the universal sheaf over $S\times \M$
(Theorem \ref{thm-diagonal}).
In section 
\ref{sec-vanishing-odd-cohomology} we 
prove Theorem \ref{thm-introduction-integral-generators}, 
finding generators for 
$H^*(\M,\Integers)$, when
$S$ is $K3$ or a rational Poisson surface. These are the Poisson surfaces 
with vanishing odd cohomology groups.
We treat these special cases 
first, because the argument is simple, yet it illustrates the general idea. 
Everything follows from the formula for the class of the diagonal
in $\M\times \M$ and the K\"{u}nneth decomposition of the 
universal sheaf in $K^0_{top}(S\times \M)$. 
In section \ref{sec-a-decomposition-of-the-diagonal-in-the-chow-ring} 
we establish the isomorphism 
$A^*(\M)\cong H^*(\M,\Integers)$, between the Chow ring and the integral 
cohomology ring, when $S$ is a rational Poisson surface
(Theorem \ref{thm-introduction-chow-ring-isomorphic-to-cohomology-ring}).
The results in sections
\ref{sec-vanishing-odd-cohomology} and 
\ref{sec-a-decomposition-of-the-diagonal-in-the-chow-ring}, 
and their proofs, are natural extensions of those of 
Ellingsrud and Str{\o}mme in the case $S=\PP^2$
\cite{ellingsrud-stromme-p2}. 
Available to us is the formula for the class of the
diagonal in $\M\times \M$, provided by Theorem
\ref{thm-diagonal}, which was proven in \cite{ellingsrud-stromme-p2}
in the special case of $\PP^2$. 
In section \ref{sec-non-simply-connected} we treat Poisson surfaces with 
non-vanishing odd cohomology groups. 

%***************************************************************************
% 
%***************************************************************************
\subsection{Stable sheaves and their moduli spaces}
\label{sec-stable-sheaves-and-modui-spaces}

Let $S$ be a smooth and projective symplectic or Poisson surface and
$H$ an ample line bundle on $S$. 
The Hilbert polynomial of a coherent sheaf $F$ on $S$ is defined by
\[
P_F(n) \ \ := \ \ \chi(F\otimes H^n) := h^0(F\otimes H^n)-
h^1(F\otimes H^n)+h^2(F\otimes H^n).
\]
Let $r$ be the rank of $F$, $f_i:=c_i(F)$, $i=1,2$, its Chern classes, 
$h:=c_1(H)$, and $K:=c_1(T^*S)$. Hirzebruch-Riemann-Roch yields
the equality
\begin{eqnarray*}
P_F(n) & = & (rh^2/2)n^2 + [(h\cdot f_1)-r/2(h\cdot K)]n +
%\\
%& & 
\left[(f_1^2-2f_2)-f_1\cdot K+2r\chi(\StructureSheaf{S})\right]/2.
\end{eqnarray*}
The degree $d$ of $P_F(n)$ is equal to 
the dimension of the support of $F$.
Let $l_0(F)/d!$ be 
the coefficient of $n^d$. 
%we get integers $l_i$, $0\leq i\leq d$, defined by the equation
%\[
%P_F(n) = \frac{l_0}{d!}n^d + \frac{l_1}{(d-1)!}n^{d-1}+ \dots + l_d(F).
%\]
Then $l_0(F)$ is a positive integer. Explicitly, 
if $r>0$ then $l_0(F) := rh^2$. If $r=0$ and $d=1$ then
$l_0(F) :=h\cdot c_1(F)$. If $r=0$ and $d=0$, then 
$l_0(F) :=-c_2(F)$.
Given two polynomials $p$ and $q$ with real coefficients,
we say that $p \succ q$ (resp. $p \succeq q$) if
$p(n) > q(n)$ (resp. $p(n) \geq q(n)$) for all $n$ sufficiently large.

%***********************************
% Definition of stability
%***********************************
\begin{defi}
\label{def-stability}
%\begin{enumerate}
%\item
%\label{def-item-gieseker-stability}
A coherent sheaf $F$ on $S$ is called {\em $H$-semi-stable} (resp. 
{\em $H$-stable}) if it has support of pure dimension 
and any non-trivial subsheaf $F' \subset F$, $F' \neq (0)$, $F' \neq F$ 
satisfies 
\[
\frac{P_{F'}}{l_0(F')} \preceq \frac{P_F}{l_0(F)} \ \ \ (\mbox{resp.} \prec).
\]
%\item
%\label{def-item-slope-stability}
%A coherent sheaf $F$ on $S$ is called {\em $H$-slope-semi-stable} 
%if it has support of pure dimension $d\geq 1$ 
%and for any non-trivial subsheaf $F'$ we have
%\[
%\frac{l_1(F')}{l_0(F')} 
%\leq
%\frac{l_1(F)}{l_0(F)}. 
%\]
%If $d=1$, $H$-slope-stability is defined using above a strict inequality.
%If $d=2$, $H$-slope-stability is defined using above a strict inequality
%and considering only non-trivial subsheaves $F'$ of lower rank. 
%\end{enumerate}
\end{defi}

Let 
$v\in K^0_{top}(S)$ be a class of rank 
$r\geq 0$ and first Chern class 
$c_1\in H^2(S,\Integers)$ of Hodge type $(1,1)$. 
The moduli space $\M_{H}(v)$, of isomorphism classes of 
$H$-stable sheaves of class $v$, 
is a quasi-projective scheme \cite{gieseker,simpson}. 

\begin{defi}
\label{def-v-suitable}
{\rm 
An ample line bundle $H$ is said to be $v$-{\em generic}, 
if every $H$-semi-stable sheaf with class $v$ is $H$-stable. 
}
\end{defi}

The class $v$ is {\em primitive}, 
if $v$ is not a multiple of another class in $K^0_{top}S$, 
by an integer larger than $1$. 
If $v$ is primitive, then 
a $v$-generic polarization exists when $r>0$, or when $r=0$ and 
$\chi(v)\neq 0$
(see \cite{yoshioka-chamber-str}, when $r>0$,  
\cite{yoshioka-abelian-surface} Lemma 1.2, when $r=0$,
and \cite{markman-reflections} Condition 3.1, 
for an existence criterion when $r=0$ and $\chi(v)=0$).
If $H$ is $v$-generic, then $\M_H(v)$ is projective
\cite{gieseker,simpson}. 
{\em Caution:} The standard definition of the term $v$-{\em generic}
is more general, and does not assume that $v$ is primitive.
We will assume, throughout the paper, the following:

\begin{condition}
\label{cond-v-suitable}
The class $v$ is primitive and $H$ is $v$-generic. 
\end{condition}

\medskip
%The connectedness assumption in Condition
%\ref{cond-v-suitable} is automatically satisfied
%for K3 and abelian surfaces \cite{yoshioka-abelian-surface} .
%If disconnected (???), 
%we may replace $\M_H(v)$ by any of its connected components.
When $S$ is a K3 or abelian surface, 
the moduli space $\M_{H}(v)$ is smooth and holomorphic symplectic
\cite{mukai-symplectic-structure}. 
It is also connected (\cite{yoshioka-abelian-surface}, or
Corollary \ref{cor-connected} below).
When $S$ is a non-symplectic Poisson surface, i.e., 
one with a non-trivial and effective anti-canonical 
line-bundle, then the moduli space $\M_{H}(v)$ is smooth, whenever the
rank $r$ is positive. When $r=0$, a moduli space $\M_{H}(v)$, 
of stable sheaves with pure one-dimensional support, may be
singular (see Example 8.6 in \cite{cime}). Smoothness 
and connectedness of $\M_H(v)$ are guaranteed, 
for sheaves with support of dimension $1$ or $2$, 
by the condition:

\begin{condition}
\label{cond-proper-subsheaf}
The sheaf $F\otimes K_S$ is isomorphic to a proper subsheaf of $F$,
for every stable sheaf $F$ parametrized by $\M_H(v)$.
\end{condition}

\medskip
The above condition and Serre's duality imply, 
that the extension group $\Ext^2(E,F)$ vanishes, for every two stable 
sheaves $E$, $F$ in $\M_H(v)$. The condition is automatically 
satisfied when $r>0$ and $S$ is a non-symplectic Poisson surface. 
The condition is also satisfied, for example, 
for sheaves with support of dimension $1$ on a del-Pezzo surface $S$
(with an ample anti-canonical line bundle $K_S^{-1}$).
Smoothness follows from Condition
\ref{cond-proper-subsheaf}, by a criterion of Artamkin \cite{artamkin}.
Connectedness is proven below (Corollary \ref{cor-connected}).

A Poisson structure on $S$ determines a Poisson structure on 
the smooth moduli space $\M_H(v)$  (see \cite{tyurin} for the
definition of the tensor, 
\cite{bottacin} for the Jacobi identity when $r>0$, and
\cite{hm-sklyanin} when $r=0$). 

%***************************************************************************
% The class of the diagonal
%***************************************************************************
\subsection{The class of the diagonal}
\label{sec-class-of-diagonal}

Assume that $S$, $v$, and $H$ satisfy Condition \ref{cond-v-suitable}.
The following result holds, when $S$ is a K3 or abelian surface
(by \cite{markman-diagonal}),
or if $S$ is a non-symplectic Poisson surface and $\M_H(v)$ satisfies 
condition
\ref{cond-proper-subsheaf} (by \cite{beauville-diagonal}).
Let $\pi_{ij}$ be the projection from 
$\M_H(v)\times S \times \M_H(v)$ onto the product of the $i$-th and $j$-th
factors. Assume, that there exists a universal sheaf over $S\times \M_H(v)$.

\begin{thm} \cite{beauville-diagonal,markman-diagonal}
\label{thm-diagonal}
Let $\E_v'$, $\E_v''$ be any two universal families of sheaves over the 
$m$-dimensional moduli space $\M_H(v)$. 
\begin{enumerate}
\item
\label{thm-item-class-of-diagonal}
The class of the diagonal, in the Chow ring of 
$\M_H(v)\times \M_H(v)$, is identified by 
\begin{equation}
\label{eq-class-of-diagonal}
c_m\left[- \ 
\pi_{13_!}\left(
\pi_{12}^*(\E'_v)^\vee\stackrel{L}{\otimes}\pi_{23}^*(\E''_v)
\right)
\right],
\end{equation}
where 
%$\pi_{13_!}$ is the K-theoretic push-forward, and 
both the dual
$(\E'_v)^\vee$ and the tensor product are taken in the derived category.
\item
When $S$ is a $K3$ or abelian surface, the following vanishing holds
\begin{equation}
\label{eq-vanishing-ofc-m-1}
c_{m-1}\left[- \ 
\pi_{13_!}\left(
\pi_{12}^*(\E'_v)^\vee\stackrel{L}{\otimes}\pi_{23}^*(\E''_v)
\right)
\right] \ \ \ = \ \ \ 0.
\end{equation}
\end{enumerate}
\end{thm}

An analogue of Theorem \ref{thm-diagonal}
holds, even when a universal sheaf does 
not exist (see Proposition 
\ref{prop-class-of-diagonal-without-universal-sheaf} below). 

The class (\ref{eq-class-of-diagonal}) is independent of the
choice of the universal sheaves $\E'_v$ and $\E''_v$, as a consequence of 
part \ref{thm-item-class-of-diagonal} of the Theorem. 
This independence is proven by another method in Lemma 
\ref{lemma-invariance-of-c-m-under-tensorization}
equation (\ref{eq-c-r-plus-1-is-invariant}), 
in the non-symplectic Poisson case. 
In the symplectic case, Lemma 
\ref{lemma-invariance-of-c-m-under-tensorization} equation
(\ref{eq-c-r-plus-2-depends-linearly-on-F}) 
relates the independence of the class (\ref{eq-class-of-diagonal}) 
to the vanishing (\ref{eq-vanishing-ofc-m-1}).

\begin{new-lemma}
\label{lemma-invariance-of-c-m-under-tensorization}
Let $X$ be a topological space, $x$ a class of rank $r\geq 0$ in 
$K^0_{top}(X)$, and $L$ a complex line-bundle on $X$. 
Then the Chern classes of $x\cup L$ satisfy the following equations
\begin{eqnarray}
\label{eq-c-r-plus-k}
c_{r+n}(x\cup L) &  = & c_{r+n}(x) 
-(n-1)c_{r+n-1}(x)c_1(L)
+ \cdots 
\\
\nonumber
& & + 
(-1)^d\Choose{n-1}{d}c_{r+n-d}(x)c_1(L)^d + \cdots +
(-1)^{n-1}c_{r+1}(x)c_1(L)^{n-1},
\end{eqnarray}
for $n\geq 1$. In particular, 
\begin{eqnarray}
\label{eq-c-r-plus-1-is-invariant}
c_{r+1}(x\cup L) & = & c_{r+1}(x),
\\
\label{eq-c-r-plus-2-depends-linearly-on-F}
c_{r+2}(x\cup L) & = & c_{r+2}(x) - c_{r+1}(x)c_1(L).
\end{eqnarray}
\end{new-lemma}
%Lemma \ref{lemma-invariance-of-c-m-under-tensorization} 
%is proven in the Appendix. It is used in the proof of
%Proposition \ref{prop-class-of-diagonal-without-universal-sheaf}.
%
%Let us observe first, that 
%equation (\ref{eq-c-r-plus-1-is-invariant}) follows from 
%the Porteous formula. 
%The Porteous formula interprets the class $c_{r+1}(x\cup L)$
%geometrically, as the class Poincare dual to a degeneracy locus of
%a regular section of the homomorphism bundle $\Hom(F,E)$, 
%where the class $x$ is $[E]-[F]$. The degeneracy locus is
%independent of $L$, since $\Hom(F\otimes L,E\otimes L)=\Hom(F,E)$. 
%When $c_{r+1}(x)=0$ and $c_{r+k}(x)\neq 0$, for some $k>1$, 
%then $\Hom(F,E)$ does not admit any regular section.
%Indeed, a regular section $\varphi$ would have constant rank and the
%class $x$ would be represented by the rank $r$ vector bundle 
%$coker(\varphi)$,
%implying the vanishing of $c_{r+k}(x)$, for all $k\geq 1$.
\noindent
{\bf Proof:}
Every element of $K_{top}^0(X)$ is of the form $[E]-[F]$, where
$E$ and $F$ are vector bundles on $X$ (see \cite{atiyah-book} section 2.1).
Let $e$ and $f$ be the ranks of $E$ and $F$, so that $r=e-f$.
Using the {\em splitting principle}, we can write the Chern polynomials
as products
$c_t(E)  =  \prod_{i=1}^e(1+\alpha_{i}t)$ and 
$c_t(F) =  \prod_{j=1}^f(1+\beta_{j}t)$,
where the $\alpha_i$ and $\beta_j$ are formal variables, and only 
the  symmetric polynomials in the $\alpha_i$ or $\beta_j$ are
interpreted as cohomology classes. Set $\ell:=c_1(L)$.
The Chern polynomial of $x\cup L$ then satisfies:
\begin{eqnarray*}
c_t(x\cup L) & = &
\prod_{i=1}^e(1+[\alpha_{i}+\ell]t)\ \ / \ \ 
\prod_{j=1}^f(1+[\beta_{j}+\ell]t)  \ \ = \ \ 
(1+\ell t)^rc_{\frac{t}{1+\ell t}}(x)  \ \ = \ \ 
\\
 & = & \sum_{q}c_q(x)t^q(1+\ell t)^{r-q}.
\end{eqnarray*}
Thus, ${\displaystyle c_{r+n}(x\cup L)=
\sum_{i=0}^{r+n}\Choose{i-n}{i}c_{r+n-i}(x)\ell^i}$. 
The lemma follows from the vanishing of $\Choose{i-n}{i}$, for $i\geq n$. 
\EndProof

Connectedness of $\M_H(v)$ is an immediate corollary of
Theorem \ref{thm-diagonal}. This was observed by Mukai, for
$2$-dimensional moduli spaces over $K3$ surfaces \cite{mukai-hodge}, 
and by Kaledin, Lehn, and Sorger,  for more general moduli spaces over
$K3$ and abelian surfaces \cite{kaledin-lehn-sorger}.

\begin{cor}
\label{cor-connected}
$\M_H(v)$ is connected.
\end{cor}

\noindent
{\bf Proof:}
Let $M$ be a connected component of $\M_H(v)$ and
$f_i$, $i=1,2$, the projection from $S\times M$ onto the $i$-th factor.
Assume first, that a universal sheaf $\E_v$ exists over
$S\times M$. 
Let $F$ be a sheaf on $S$ with class $v$ and 
set
\[
x \ \ \ := \ \ \ 
- \ 
f_{2_!}\left[
f_{1}^*(F)^\vee\stackrel{L}{\otimes}f_{2}^*(\E_v)
\right].
\]
The class $c_m(x)$ depends only on the class $v$ and is independent 
of the sheaf $F$. 
If $F$ belongs to $M$, then $c_m(x)$ is 
Poincare dual to the pullback of the class
(\ref{eq-class-of-diagonal}), via the embedding
of $M$ in $\M_H(v)\times \M_H(v)$
sending $E$ to $(F,E)$. Thus  $c_m(x)$ is Poincare-dual to a point.

Assume that there exists an $H$-stable sheaf $F$, which 
does not belong to $M$. Then the higher-direct images 
$R^{i}f_{2_*}[f_{1}^*(F)^\vee\stackrel{L}{\otimes}f_{2}^*(\E_v)]$ 
vanish, for $i=0,2$. Consequently, 
the class $x$ is represented by the locally free sheaf
$R^{1}f_{2_*}[f_{1}^*(F)^\vee\stackrel{L}{\otimes}f_{2}^*(\E_v)]$ of
rank $-\chi(v^\vee\otimes v)$. This rank is 
$m-2$, if $S$ is symplectic, and $m-1$, if $S$ is non-symplectic but Poisson.
Thus $c_m(x)$ vanishes. 
This contradicts the non-vanishing of $c_m(x)$ proven above, 
so such $F$ can not exist.

A  universal sheaf $\E_v$ 
exists always over $S\times \PP$, where $\PP$ is 
a projective bundle over $\M_H(v)$
%even if a universal sheaf does not exist over $\M_H(v)$ 
given in equation (\ref{eq-PP}). 
The above argument generalizes, as it proves that
$c_m(x)=0$, when $F\not\in M$, 
and $c_m(x)$ is Poincare-dual to a fiber of $\PP$
over a point of $M$, when $F\in M$. 
\EndProof

%**************************************************************
% Poisson surfaces with vanishing odd cohomology
%**************************************************************
\subsection{Poisson surfaces with vanishing odd cohomology}
\label{sec-vanishing-odd-cohomology}
Assume that a universal sheaf $\E$ exists over $S\times \M_H(v)$.
The assumption is dropped in section \ref{sec-twisted-sheaves}. 
Assume, in addition, the following:

\begin{condition}
\label{cond-odd-cohomology-vanishes}
The cohomology groups $H^i(S,\Integers)$ vanish, for odd $i$.
\end{condition}

\noindent
Equivalently, $S$ is either a K3, 
or a smooth projective rational Poisson surface. 
The case of a general Poisson surface will be treated in section
\ref{sec-non-simply-connected}.
Condition \ref{cond-odd-cohomology-vanishes} implies, that 
$K^1_{top}(S)$ vanishes, the group 
$K^0_{top}(S)$ is free of rank equal to
$H^*(S,\Integers)$, and the Chern character
\[
ch \ : \ K^0_{top}(S) \ \ \longrightarrow \ \ 
H^*(S,\RationalNumbers)
\]
is an injective homomorphism 
(\cite{atiyah-hirzebruch} page 19).
%When $S$ is a $K3$ surface, $ch$ induces an isomorphism
%$K^0_{top}(S)\cong H^*(S,\Integers)$ (the Chern character 
%has values in $H^*(S,\Integers)$, since $H^2(S,\Integers)$ is an
%{\em even} lattice and surjectivity of $ch$ is easily verified). 
%This special property of 
%$K3$ surfaces will be used only in (???).

Given any cell complex $M$, the K\"{u}nneth Theorem provides an isomorphism
\[
K^0_{top}(S)\otimes K^0_{top}(M) \ \ \cong \ \ 
K^0_{top}(S\times M),
\]
given by the exterior product 
(\cite{atiyah-book}, Corollary 2.7.15). We used here the vanishing of 
$K^1_{top}S$. Use a basis  $\{x_1, \dots, x_n\}$ of $K^0_{top}(S)$
to write the K\"{u}nneth decomposition of the class of the universal sheaf:
\begin{equation}
\label{eq-e-i}
\E \ \ \equiv \ \ 
\sum_{i=1}^n x_i\otimes e_i.
\end{equation}

\begin{prop}
\label{prop-integral-generators}
\begin{enumerate}
\item
\label{cor-item-integral-generators}
The cohomology ring $H^*(\M_H(v),\Integers)$ is generated by the
Chern classes $c_j(e_i)$ of the classes $e_i\in K^0_{top}\M_H(v)$, 
which are given in equation (\ref{eq-e-i}).
\item
\label{cor-item-odd-cohomology-vanishes}
The cohomology group $H^i(\M_H(v),\Integers)$ 
vanishes, for odd $i$, and is torsion free, when $i$ is even.
\end{enumerate}
\end{prop}

\noindent
{\bf Proof:}
\ref{cor-item-integral-generators}) 
Set $\M:=\M_H(v)$. 
The proof is similar to that of Theorem 2.1 in \cite{ellingsrud-stromme-p2}. 
The projection formula (\ref{eq-projection-formula})
yields the equivalence in 
$K^0_{top}(\M\times \M)$
\begin{equation}
\label{eq-first-decomposition-of-diagonal}
- \ 
\pi_{13_!}\left(
\pi_{12}^*(\E)^\vee\stackrel{L}{\otimes}\pi_{23}^*(\E)
\right) \ \ \equiv \ \ 
-\sum_{i=1}^n \sum_{j=1}^n p_1^!(e_i^\vee)\cup p_2^!(e_j)\cup 
\pi_{13_!}\pi_2^!(x_i^\vee\cup x_j),
\end{equation}
where $p_i$ is the projection from $\M\times \M$ on
the $i$-th factor, $i=1,2$, and $\pi_{ij}$ are the projections
from $\M\times S\times \M$. 
Property 
(\ref{eq-cohomology-and-base-change}) of the Gysin homomorphism 
implies the equality 
$\pi_{13_!}\pi_2^!(x_i^\vee\cup x_j)=\chi(x_i^\vee\cup x_j)\cdot 1$. 
We define the {\em Mukai pairing} on $K^0_{top}(S)$ by
\begin{equation}
\label{eq-Mukai-pairing}
(x,y) \ \ := \ \ -\chi(x^\vee\cup y). 
\end{equation}
It is known to be a perfect pairing (see also 
Remark \ref{rem-perfect-pairing}).
The equivalence (\ref{eq-first-decomposition-of-diagonal}) becomes
\[
- \ 
\pi_{13_!}\left(
\pi_{12}^*(\E)^\vee\stackrel{L}{\otimes}\pi_{23}^*(\E)
\right) \ \ \equiv \ \ 
\sum_{i=1}^n \sum_{j=1}^n 
(x_i,x_j)
p_1^!(e_i^\vee)\cup p_2^!(e_j).
\]
Theorem \ref{thm-diagonal} translates to the equality
\begin{equation}
\label{eq-translation-of-thm-diagonal}
\delta \ \ = \ \ c_m\left(
\sum_{i=1}^n  \sum_{j=1}^n 
(x_i,x_j) p_1^!(e_i^\vee)\cup p_2^!(e_j)
\right),
\end{equation}
where $\delta$ is Poincare-dual to the class of the diagonal 
in $\M\times \M$.
Clearly, the $m$-th Chern class, on the right hand side,
can be written as a sum 
\begin{equation}
\label{eq-decomposition-of-the-diagonal-of-moduli}
\delta \ \ = \ \ \sum_{j\in J}p_1^*\alpha_j \cup p_2^*\beta_j,
\end{equation}
where each $\alpha_j$ and $\beta_j$ is a polynomial, with
integral coefficients,  
in the Chern classes of $e_i$, for $1\leq i \leq n$. 

We have the well known formula
\begin{equation}
\label{eq-gysin-composed-with-cup-delta-is-identity}
x \ \ = \ \ p_{1,*}(\delta\cup p_2^*x), 
\end{equation}
where $x$ is a class in $H^*(\M,\Integers)$.  
We recall the proof of that formula. 
Let $[\M]\in H_{2m}(\M,\Integers)$ and 
$[\M\times \M]\in H_{4m}(\M\times \M,\Integers)$ be the orientation classes 
and $[\Delta]\in H_{2m}(\M\times \M,\Integers)$ the class of the diagonal. 
The Poincare Duality isomorphism, sending $\delta$ to $[\Delta]$, 
is induced by the cap product
\[
(\bullet) \cap [\M\times \M] \ : \ H^*(\M\times \M,\Integers) \ \ 
\LongIsomRightArrow \ \ H_*(\M\times \M,\Integers).
\]
The Poincare dual of equation 
(\ref{eq-gysin-composed-with-cup-delta-is-identity}) follows
from the following equalities.
\begin{eqnarray*}
p_{1,*}\left\{(\delta\cup p_2^*x)\cap [\M\times \M]\right\}
& = &
p_{1,*}\left\{p_2^*x\cap (\delta\cap [\M\times \M])\right\}
\ = \
p_{1,*}\left\{p_2^*x\cap [\Delta]\right\}
\ = \
\\
& = & 
p_{1,*}\left\{p_1^*x\cap [\Delta]\right\}
\ = \
x\cap [\M].
\end{eqnarray*}

Using equation (\ref{eq-gysin-composed-with-cup-delta-is-identity}), 
the projection formula, 
and the decomposition 
(\ref{eq-decomposition-of-the-diagonal-of-moduli}),
%$\delta=\sum_{j\in J}p_1^*\alpha_j p_2^*\beta_j$, 
we express 
$x$ as a linear combination of the $\alpha_j$'s:
\begin{equation}
\label{eq-x-is-a-linear-combination-of-alpha-j}
x \ \ = \ \ \sum_{j\in J} \left(\int_\M x\cup\beta_j\right)\alpha_j.
\end{equation}
\ref{cor-item-odd-cohomology-vanishes})
Part \ref{cor-item-odd-cohomology-vanishes} follows from the vanishing of
the coefficient $\int_\M x\cup\beta_j$ in 
(\ref{eq-x-is-a-linear-combination-of-alpha-j}), 
whenever $x$ is a torsion class or a class of odd degree.
\EndProof

\begin{example}
{\rm
Let $S=\PP^2$ with the basis $\{y_1,y_2,y_3\}$ of $K^0_{alg}\PP^2$ given by 
$y_i:=\StructureSheaf{\PP^2}(-i)$. 
The dual basis $\{x_1,x_2,x_3\}$, with respect to 
the pairing $\chi(x_i\cup y_j)$, is given by 
$x_1:=\StructureSheaf{\PP^2}(-2)-3\StructureSheaf{\PP^2}(-1)-
\StructureSheaf{\PP^2}$, 
$x_2:=\StructureSheaf{\PP^2}(-1)-3\StructureSheaf{\PP^2}$, and 
$x_3:=\StructureSheaf{\PP^2}$. 
It is well-known, that $\alpha:K^0_{alg}\PP^2\rightarrow K^0_{top}\PP^2$ 
is an isomorphism 
(see part \ref{lemma-item-alg-and-top-K-groups-are-isomorphic} 
of Lemma \ref{lemma-algebraic-kunneth-decomposition}). 
Let $f_i$ be the projection from $\PP^2\times \M$
onto the $i$-th factor, $i=1,2$. 
Then $e_i$ in (\ref{eq-e-i}) is equivalent to
$f_{2,!}(\E\cup f_1^!y_i)$. 
We conclude, that the classes $e_i$ in Proposition
\ref{prop-integral-generators} are given by 
$e_i=f_{2,!}(\E\otimes f_1^!\StructureSheaf{\PP^2}(-i))$, 
for $1\leq i\leq 3$.
These are precisely the classes in $K^0_{alg}\M$ chosen by 
Ellingsrud and Str{\o}mme in their version of the statement of part 
\ref{cor-item-integral-generators} of Proposition
\ref{prop-integral-generators} above 
(they worked with $K^0_{alg}$ and the Chow-ring,
instead of $K^0_{top}$ and integral cohomology). 
The similarity is not a coincidence; the proof of Proposition
\ref{prop-integral-generators} is modeled after that of
Theorem 2.1 in \cite{ellingsrud-stromme-p2}.
}
\end{example}

%*************************************************************
% A Chow ring decomposition of the diagonal 
%*************************************************************
\subsection{A decomposition of the diagonal in the Chow ring}
\label{sec-a-decomposition-of-the-diagonal-in-the-chow-ring}
We prove Theorem
\ref{thm-introduction-chow-ring-isomorphic-to-cohomology-ring} in this section.
Assume, that $S$ is a smooth, projective, rational and 
Poisson surface. 
%In particular, $S$ could be any minimal rational surface.
%with vanishing odd cohomology groups $H^i(S,\Integers)$, $i=1,3$. 
Let $v\in H^*(S,\Integers)$ be a vector satisfying condition
\ref{cond-proper-subsheaf}. 
The cohomological decomposition 
(\ref{eq-decomposition-of-the-diagonal-of-moduli}), of the diagonal class, 
has an analogue in the Chow ring  $A^*(\M\times \M)$ of $\M\times \M$ 
%whenever the K\"{u}nneth formula holds for $K^0_{alg}(S\times S)$
(Theorem
\ref{thm-introduction-chow-ring-isomorphic-to-cohomology-ring}). 

\begin{defi}
\label{def-decomposition-of-diagonal}
{\rm
Let $X$ be an algebraic variety and $\Delta\subset X\times X$
the diagonal. 
We say, that  $X$ 
{\em 
admits a decomposition of the diagonal in $K^0_{alg}(X\times X)$,
}
if the class of $\StructureSheaf{\Delta}$ is in the image of
$K^0_{alg}X\otimes K^0_{alg}X\rightarrow K^0_{alg}(X\times X)$.
%We say that $X$ 
%{\em 
%admits a decomposition of the diagonal in $A^*(X\times X)$,
%}
%if the class of the diagonal is in the image of
%$A^*(X)\otimes A^*(X) \rightarrow A^*(X\times X)$.
}
\end{defi}

The decomposition in the definition exists, when $X$ is a projective space 
$\PP^n$ or the projectivization $\PP{E}$ of a vector bundle $E$ over $\PP^n$
(\cite{fulton} Example 15.1.1). Special cases include the
minimal rational surfaces: $\PP^2$, $\PP^1\times \PP^1$, and the
Hirzebruch surfaces.

\begin{new-lemma}
\label{lemma-decomposition-of-the-diagonal-for-rational-surfaces}
Every smooth and projective rational surface 
$S$ admits a decomposition of the diagonal in 
$K^0_{alg}(S\times S)$.
\end{new-lemma}

%\noindent
%{\bf Proof:}
%$S$ is obtained from $\PP^2$, $\PP^1\times \PP^1$, or a
%Hirzebruch surface, by a sequence of blow-ups. 
%\EndProof
Lemma \ref{lemma-decomposition-of-the-diagonal-for-rational-surfaces}
follows from Lemma \ref{lemma-decomposition-of-the-diagonal-after-blow-up}.
Let $X$ be a smooth complex algebraic surface, 
$\beta:\widehat{X}\rightarrow X$ the blow-up of $X$ at a closed point,
$E\subset \widehat{X}$ the exceptional divisor, 
$\delta:\Delta\rightarrow X\times X$, 
$\hat{\delta}:\widehat{\Delta}\rightarrow \widehat{X}\times \widehat{X}$, 
the diagonals, and 
$\pi_i:\widehat{X}\times \widehat{X}\rightarrow \widehat{X}$, $i=1,2$, 
the projections.

\begin{new-lemma}
\label{lemma-decomposition-of-the-diagonal-after-blow-up}
The following equivalence holds in $K^0_{alg}(\widehat{X}\times \widehat{X})$.
\begin{equation}
\label{eq-decomposition-of-diagonal-in-blow-up}
\hat{\delta}_!\StructureSheaf{\widehat{\Delta}} \ \ \  \equiv \ \ \ 
(\beta\times\beta)^![\delta_!\StructureSheaf{\Delta}]
- \left\{
\pi_1^!\left[\StructureSheaf{\widehat{X}}(E)-
\StructureSheaf{\widehat{X}}\right]
\cup \pi_2^!\left[\StructureSheaf{\widehat{X}}(E)-
\StructureSheaf{\widehat{X}}\right]
\right\}.
\end{equation}
Consequently, 
if $X$ admits a decomposition of the diagonal 
in $K^0_{alg}(X\times X)$, then the same is true 
when we replace $X$ by its blow-up $\widehat{X}$.
\end{new-lemma}

Note, that $(\beta\times \beta)_!\circ (\beta\times \beta)^!=id$,
$\eta:=(\beta\times \beta)^!\circ (\beta\times \beta)_!$ is an idempotent, 
and the right hand side of 
(\ref{eq-decomposition-of-diagonal-in-blow-up}) 
is compatible with the decomposition ${\mbox Im}(\eta)\oplus \ker(\eta)$
of $K^0_{alg}(\widehat{X}\times \widehat{X})$. 

\noindent
{\bf Proof of Lemma \ref{lemma-decomposition-of-the-diagonal-after-blow-up}:}
Let $\Delta_E\subset E\times E$ be the the diagonal.
Set $Z:=\widehat{X}\cup[E\times E]$. 
The inclusions of $E\times E$, $Z$, and $\Delta_E$
into $\widehat{X}\times \widehat{X}$ are all denoted by $\iota$. 
The desired decomposition (\ref{eq-decomposition-of-diagonal-in-blow-up})
follows from the following equivalences:
\begin{eqnarray}
\label{eq-vanishing-of-higher-torsion-sheaves}
\iota_!\StructureSheaf{Z} & \equiv &
(\beta\times\beta)^![\delta_!\StructureSheaf{\Delta}].
\\
\label{eq-inclusion-exclusion}
\iota_!\StructureSheaf{Z} & \equiv & 
\hat{\delta}_!\StructureSheaf{\widehat{\Delta}}+ 
\iota_!\StructureSheaf{E\times E} - \iota_!\StructureSheaf{\Delta_E}.
%\\
%\hat{\delta}_!\StructureSheaf{\widehat{\Delta}} & \equiv &
%(\beta\times\beta)^![\delta_!\StructureSheaf{\Delta}]
%-\iota_!\StructureSheaf{E\times E} + \iota_!\StructureSheaf{\Delta_E}.
\\
\label{eq-decomposition-of-E-times-E}
\iota_!\StructureSheaf{E\times E} & \equiv & 
\pi_1^!\left[1-\StructureSheaf{\widehat{X}}(-E)\right]\cup 
\pi_2^!\left[1-\StructureSheaf{\widehat{X}}(-E)\right].
\\
\label{eq-decomposition-of-diagonal-of-E}
\iota_!\StructureSheaf{\Delta_E}  & \equiv &
\left\{
1-[\pi_1^!\StructureSheaf{\widehat{X}}(E)\cup
\pi_2^!\StructureSheaf{\widehat{X}}(E)]
\right\}\cup \iota_!\StructureSheaf{E\times E}.
\end{eqnarray}

We prove 
(\ref{eq-vanishing-of-higher-torsion-sheaves}) first.
$Z$ is the total transform of $\Delta$ in $\widehat{X}\times \widehat{X}$
and $\iota_*\StructureSheaf{Z}=
(\beta\times\beta)^*(\delta_*\StructureSheaf{\Delta})$.
We need to prove the vanishing of the torsion sheaves
$\SheafTor_i^{\StructureSheaf{X\times X}}(
\StructureSheaf{\widehat{X}\times \widehat{X}},
\delta_*\StructureSheaf{\Delta})$, $i\geq 1$. 
Both $\Delta$ and $Z$ have codimension $2$. $\Delta$ is smooth
and is hence a local complete intersection.
Hence, so is $Z$. Furthermore, 
local equations $g_1$, $g_2$ of $\Delta$ in $X\times X$ pull back 
to local equations of $Z$, which form a regular sequence
(\cite{matsumura} Ch. 6 section 16 Theorem 31). Thus,
the Koszul complex, 
locally resolving the sheaf $\delta_*\StructureSheaf{\Delta}$, 
pulls back via $\beta\times \beta$ to a Koszul complex 
locally resolving the sheaf 
$\iota_*\StructureSheaf{Z}$ 
(\cite{matsumura} Ch. 7 section 18 Theorem 43). The vanishing of
the higher torsion sheaves follows.

Let $I_{\Delta_E}$ be the ideal sheaf of $\Delta_E$ in $E\times E$.
The embedding $E\times E\hookrightarrow Z$ pushes forward $I_{\Delta_E}$
to the ideal sheaf of $\widehat{\Delta}$ in $Z$.
Equivalence (\ref{eq-inclusion-exclusion}) follows. 
Equivalence (\ref{eq-decomposition-of-E-times-E}) is clear.
Note, that $\StructureSheaf{\widehat{X}}(E)$ restricts to $E$ as
$\StructureSheaf{E}(-1)$. Hence, 
$1-[\pi_1^!\StructureSheaf{\widehat{X}}(E)\cup
\pi_2^!\StructureSheaf{\widehat{X}}(E)]$ is sent to the class of
$\StructureSheaf{\Delta_E}$ via the composition
$K_{alg}(\widehat{X})\otimes K_{alg}(\widehat{X})\rightarrow 
K_{alg}(E)\otimes K_{alg}(E)\rightarrow K_{alg}(E\times E).$
Equivalence (\ref{eq-decomposition-of-diagonal-of-E}) now follows 
from the projection formula (\ref{eq-projection-formula}).
\EndProof

\begin{new-lemma}
\label{lemma-algebraic-kunneth-decomposition}
Let $X$ be a projective algebraic variety, which  
admits a decomposition 
\begin{equation}
\label{eq-decomposition-of-O-S}
\StructureSheaf{\Delta_X} \equiv\sum_{i\in I}p_1^!x_i\cup p_2^!y_i
\end{equation}
of the diagonal in $K^0_{alg}(X\times X)$. 
Then the following statements hold.
\begin{enumerate}
\item
\label{lemma-item-x-i-generate-K-S}
The $x_i$ generate $K^0_{alg}X$. 
\item
\label{lemma-item-K-S-is-torsion-free}
$K^0_{alg}X$ is a finitely generated free $\Integers$-module.
\item
\label{lemma-item-dual-basis}
Suppose that the set $\{x_i\}_{i\in I}$ is minimal.
Then $\{x_i\}_{i\in I}$ and $\{y_i\}_{i\in I}$ are dual bases,
with respect to the pairing $\chi(x_i\cup y_j)$. 
\item
\label{lemma-item-kunneth-decomposition}
Exterior product
\[
\mu \ : \ K^0_{alg}(X)\otimes K^0_{alg}(M) \ \ \longrightarrow \ \ 
K^0_{alg}(X\times M)
\]
is an isomorphism, for every algebraic variety $M$. 
\item
\label{lemma-item-alg-and-top-K-groups-are-isomorphic}
$\alpha : K^0_{alg}X\rightarrow K^0_{top}X$ is an isomorphism.
\item
\label{lemma-item-same-for-chow-ring}
If $X$ is smooth and projective, and we replace assumption 
(\ref{eq-decomposition-of-O-S})
by its Chow ring analogue
\begin{equation}
\label{eq-decomposition-in-chow-ring}
[\Delta] \ \ \equiv \ \ \sum_{i\in I}p_1^*\alpha_i\cup p_2^*\beta_i
\end{equation}
in $A^*(X\times X)$, then the statements above hold, after replacing
$K^0_{alg}$ by the Chow ring, $K^0_{top}X$ by $H^*(X,\Integers)$,
and the pairing $\chi(x_i\cup y_j)$ by the intersection pairing
$\int_X \alpha_i\cup\beta_j$.
\end{enumerate}
\end{new-lemma}

\noindent
{\bf Proof:}
The proof is again similar to that of Theorem
2.1 in \cite{ellingsrud-stromme-p2}.
Let $p_i$ be the projection from $X\times X$, 
$f_i$ the projection from $X\times M$, and
$\pi_{ij}$ the projection from $X\times X \times M$ onto the product of the
$i$-th and $j$-th factors.  
Everything follows from the evident formula
\[
x \ \ \equiv \ \ p_{1,!}\left(\StructureSheaf{\Delta_X}\cup p_2^!x\right),
\]
for any $x$ in $K^0_{alg}X$. 
Using the projection formula (\ref{eq-projection-formula}) and the 
decomposition (\ref{eq-decomposition-of-O-S}), we get 
\begin{equation}
\label{eq-x-is-a-linear-combination-of-x-i}
x \ \ \equiv \ \ \sum_{i\in I}\chi(x\cup y_i)x_i.
\end{equation}

\ref{lemma-item-x-i-generate-K-S}) 
Part \ref{lemma-item-x-i-generate-K-S} follows from equation 
(\ref{eq-x-is-a-linear-combination-of-x-i}).

\ref{lemma-item-K-S-is-torsion-free})
For part \ref{lemma-item-K-S-is-torsion-free} observe, that if 
$x\in K^0_{alg}X$
is a torsion class, then the coefficients 
$\chi(x\cup y_i)$ in equation (\ref{eq-x-is-a-linear-combination-of-x-i}) 
vanish for all $i$. 

\ref{lemma-item-dual-basis})
The minimality assumption implies, that the $x_i$ are linearly independent. 
The statement follows by setting $x=x_j$ in equation 
(\ref{eq-x-is-a-linear-combination-of-x-i}). 

\ref{lemma-item-kunneth-decomposition}) 
Let $\E$ be a class in $K^0_{alg}(X\times M)$.
The projection formula (\ref{eq-projection-formula}) yields the decomposition
\[
\E \ \ \equiv \ \ 
\pi_{13_!}\left(
\pi_{12}^!\StructureSheaf{\Delta_X}\cup \pi_{23}^!\E
\right) \ \ \equiv \ \ 
\sum_{i\in I}f_1^!x_i\cup \pi_{13_!}(\pi_2^!y_i\cup \pi_{23}^!\E).
\]
The Cohomology and Base Change Theorem, applied to the fiber product
\[
\begin{array}{ccc}
X\times X\times M & \LongRightArrowOf{\pi_{23}} & X\times M
\\
\pi_{13} \ \downarrow \ \hspace{3ex} & & 
\hspace{2ex} \ \downarrow \ f_2
\\
X\times M & \LongRightArrowOf{f_2} & M,
\end{array}
\]
 implies the second equality below
\[
\pi_{13_!}(\pi_2^!y_i\cup \pi_{23}^!\E) \ \ \equiv \ \ 
\pi_{13_!}\pi_{23}^!(f_1^!y_i\cup\E) \ \ \equiv \ \ 
f_2^!f_{2_!}(f_1^!y_i\cup\E).
\]
Consequently, 
$\E\equiv \sum_{i\in I}f_1^!x_i\cup f_2^!f_{2_!}(f_1^!y_i\cup\E)$ and 
the exterior product is surjective. 

We prove injectivity next. 
Choose dual bases as in part \ref{lemma-item-dual-basis}.
Suppose $\E := \sum_{i\in I}x_i\otimes e_i$ is in the kernel of $\mu$. 
%We may assume, that the $z_i$ are linearly independent classes in
%$K^0_{alg}X$. 
The projection formula  (\ref{eq-projection-formula}) yields
\[
0 \ \ = \ \ f_{2_!}\left(f_1^!y_j\cup \mu(\E)\right)
\ \ = \ \ 
f_{2_!}\left(f_1^!y_j\cup
\left[\sum_{i\in I}f_1^!x_i\cup f_2^!e_i
\right]
\right)
\ \ = \ \ 
\sum_{i\in I}\chi(y_j\cup x_i)e_i \ \ =  \ \ e_j.
\]
Hence, all the $e_j$ vanish. 

\ref{lemma-item-alg-and-top-K-groups-are-isomorphic})
Surjectivity of 
$\alpha : K^0_{alg}X\rightarrow K^0_{top}X$ follows from equation
(\ref{eq-x-is-a-linear-combination-of-x-i}), interpreted in $K^0_{top}X$.
Injectivity follows from the vanishing of $\chi(x\cup y_i)$
in equation (\ref{eq-x-is-a-linear-combination-of-x-i}), for a class $x$
in the kernel of $\alpha$.

\ref{lemma-item-same-for-chow-ring}) 
The proof of part \ref{lemma-item-same-for-chow-ring} is 
a straightforward translation of the above proofs. We include only the 
translation of the proof
of part \ref{lemma-item-alg-and-top-K-groups-are-isomorphic}. 
Once we replace equation (\ref{eq-decomposition-of-the-diagonal-of-moduli}) 
by equation (\ref{eq-decomposition-in-chow-ring}), then 
equation (\ref{eq-x-is-a-linear-combination-of-alpha-j})
expresses every cohomology class $x$ in $H^*(X,\Integers)$ 
as a linear combination of classes coming from $A^*(X)$. The surjectivity of
$A^*(X)\rightarrow H^*(X,\Integers)$ \ follows. 
When the class $x$ belongs to $A^*(X)$, equation 
(\ref{eq-x-is-a-linear-combination-of-alpha-j}) holds in $A^*(X)$.
The injectivity follows
from the vanishing of the coefficients $\int_{X}x\cup \beta_j$ in 
equation (\ref{eq-x-is-a-linear-combination-of-alpha-j}), when 
$x$ is a class in the kernel of $A^*(X)\rightarrow H^*(X,\Integers)$. 
\EndProof

\begin{rem}
\label{rem-perfect-pairing}
{\rm
The topological Mukai pairing (\ref{eq-Mukai-pairing}) is a perfect pairing. 
This follows from the  
argument in the proof of part \ref{lemma-item-dual-basis} of 
Lemma \ref{lemma-algebraic-kunneth-decomposition},
provided equation (\ref{eq-decomposition-of-O-S}) 
is taken in $K^0_{top}[S\times S]$. 
}
\end{rem}

\noindent
{\bf Proof of Theorem
\ref{thm-introduction-chow-ring-isomorphic-to-cohomology-ring}:}
There exists a class $x$ in $K_{alg}^0(S)$, such that
$\chi(x\cup v)=1$, by part 
\ref{lemma-item-dual-basis} of Lemma 
\ref{lemma-algebraic-kunneth-decomposition}. Hence, 
a universal sheaf exists \cite{mukai-hodge}.
The decomposition (\ref{eq-e-i}), of the universal sheaf, 
can be taken in $K^0_{alg}(S\times \M)$, by part
\ref{lemma-item-kunneth-decomposition} of Lemma 
\ref{lemma-algebraic-kunneth-decomposition}. 
Consequently, the decomposition 
(\ref{eq-decomposition-of-the-diagonal-of-moduli}), of the diagonal
in $\M\times \M$, is given in the Chow ring $A^*(\M\times \M)$.
The fact, that $A^*(\M)\rightarrow H^*(\M,\Integers)$ is an isomorphism, 
follows from 
parts \ref{lemma-item-alg-and-top-K-groups-are-isomorphic} and
\ref{lemma-item-same-for-chow-ring} of Lemma
\ref{lemma-algebraic-kunneth-decomposition}.
\EndProof

%***************************************************************************
% Non-simply connected Poisson surfaces
%***************************************************************************
\subsection{Non-simply connected Poisson surfaces}
\label{sec-non-simply-connected}

We omit Condition \ref{cond-odd-cohomology-vanishes} and consider any 
smooth, projective, symplectic or Poisson surface $S$. 
%replacing it with:
%
%\medskip
%\noindent
%{\bf Condition \ref{cond-odd-cohomology-vanishes}':}
%The cohomology groups $H^i(S,\Integers)$ are torsion free, for all $i$. 

\medskip
We need to review first some background from Topological K-Theory. 
Let $X$ be a connected topological space. 
The reduced $K$-group $\widetilde{K}^0_{top}(X)$ is the kernel of 
the restriction homomorphism $K^0_{top}(X)\rightarrow K^0_{top}(x_0)$,
where $x_0$ is a point of $X$. 
Given two topological spaces $X$ and $Y$, with points $x_0\in X$ and 
$y_0\in Y$, we set $X\vee Y:= X\times \{y_0\}\cup \{x_0\}\times Y$ and
$X\wedge Y:=X\times Y/X\vee Y$. 
Denote by $SX$ the reduced suspension of $X$,
$SX:=S^1\wedge X$,
% $SX:=S^1\times X/S^1\times\{x_0\}\cup\{s_0\}\times X$, 
where $s_0$ is a fixed point of the circle $S^1$.
Let $X^+$ be the disjoint union of
$X$ and a point $x_0$. The notation $SX^+$ stands for $S(X^+)$. 
Note, that $S(X^+)=S^1\times X/\{s_0\}\times X$.
The associativity of the operation $\wedge$ yields 
$S(S(X^+))=S^2\times X/\{s_0\}\times X$ and similarly for the 
$n$-th iterate $S^n(X^+)=S^n\times X/\{s_0\}\times X$.
Recall, that the odd $K$-group $K_{top}^1(X)$ 
is defined to be $\widetilde{K}^0_{top}(SX^+)$
(see \cite{atiyah-book}). 
The latter group is naturally isomorphic to the 
kernel of
\[
K^0_{top}(S^1\times X) \ \ \LongRightArrowOf{\iota^*} \ \ 
K^0_{top}(\{s_0\}\times X),
\]
where $s_0$ is a fixed point in $S^1$ and 
$\iota:\{s_0\}\times X\hookrightarrow S^1\times X$ is the natural embedding 
(\cite{atiyah-book}, Corollary 2.4.7).
Set $K_{top}^*(X):=K_{top}^0(X)\oplus K_{top}^1(X)$. 
Let $X$ and $Y$ be smooth complex algebraic varieties.
A proper morphisms $f:X\rightarrow Y$ extends to a proper continuous map 
$\tilde{f}:SX\rightarrow SY$. Gysin maps are defined, more generally, for 
proper morphism between differentiable manifolds, with even dimensional 
fibers, satisfying an additional condition; existence and choice of a 
certain relative $^c$spinorial structure 
(see \cite{karoubi} Proposition IV.5.24 and remark IV.5.27). 
When $f$ is a proper morphism between complex manifolds $X$ and $Y$, 
a natural $^c$spinorial structure exists for $f$, 
as well as for the suspension $\tilde{f}$ 
(use \cite{karoubi} Theorem II.4.8 for the latter). 
Consequently, we get a Gysin map 
$f_!:=\tilde{f}_!: K_{top}^1(X) \rightarrow K_{top}^1(Y)$. 
The projection formula  (\ref{eq-projection-formula})
and property (\ref{eq-cohomology-and-base-change})
extend for classes in $K_{top}^1$. 

The exterior product 
$K^0_{top}(X)\otimes K^1_{top}(Y)\rightarrow K^1(X\times Y)$
is defined as follows. 
There is a natural map $q:X\times SY^+ \rightarrow S(X\times Y)^+$. 
The image of the exterior product 
$K^0_{top}(X)\otimes \widetilde{K}^0_{top}(SY^+)\rightarrow 
\widetilde{K}^0(X\times SY^+)$
is contained in the image of $\widetilde{K}^0[S(X\times Y)^+]$ 
via $q^!$. Since $q^!$ is injective, the exterior product has values in 
$K^1(X\times Y)$. 
The exterior product 
\[
K^1_{top}(X)\otimes K^1_{top}(Y) \ \ \rightarrow \ \ K^0(X\times Y)
\]
is defined, using the natural map
$q^2:SX^+\times SY^+ \rightarrow S^2(X\times Y)^+$ and 
Bott's Periodicity Theorem,
which implies the isomorphism 
$\beta:K^0[X\times Y]\IsomRightArrow \widetilde{K}^0[S^2(X\times Y)^+]$
(see Theorem 2.4.9 and section 2.6 in \cite{atiyah-book}). 
Under the inclusion of
$\widetilde{K}^0[S^2(X\times Y)^+]$ in 
$K^0[S^2\times X\times Y]$, the isomorphism $\beta$ sends a class $\alpha$ 
in $K^0[X\times Y]$ to the exterior product of $\alpha$ with
the generator of $\widetilde{K}^0(S^2)$. 

\begin{defi}
\label{def-odd-chern-classes}
Let $x$ be a class in $K_{top}^1(X)$, corresponding to a class $\tilde{x}$ in
$\widetilde{K}^0_{top}(SX^+)$, and $i\geq 1/2$ a half-integer. 
The Chern class $c_i(x)$ of $x$ is defined as the image in 
$H^{2i}(X,\Integers)$ of $c_{i+\frac{1}{2}}(\widetilde{x})$, 
via the isomorphism 
$H^{2i}(X,\Integers)\cong H^{2i}(X^+,\Integers)\cong 
H^{2i+1}(SX^+,\Integers)$.
\end{defi}

One extends, similarly, the Chern character to $K^*_{top}X$. 
The Chern character is a linear homomorphism, which preserves products
and commutes with pullbacks. 
It maps $K^1_{top}X$ into $H^{odd}(X,\RationalNumbers)$. 
When $H^*(X,\Integers)$ is torsion free, then so is $K^*_{top}X$.
Furthermore,
\[
ch \ : \ K^*_{top}X\otimes \RationalNumbers \ \ \longrightarrow \ \ 
H^*(X,\RationalNumbers) 
\]
is an isomorphism (\cite{atiyah-hirzebruch} and \cite{karoubi} V.3.26).
The K\"{u}nneth Theorem holds as well (\cite{atiyah-book} Theorem 2.7.15). 
When $H^*(X,\Integers)$ is torsion free, it states that the exterior product
\[
K^0_{top}X\otimes K^0_{top}M \ \ \ \oplus \ \ \ 
K^1_{top}X\otimes K^1_{top}M \ \ \ 
\LongIsomRightArrow \ \ \ 
K^0_{top}[X\times M]
\]
is an isomorphism, for any cell complex $M$. 

The cohomology group $H^*(S,\Integers)$ is torsion free, 
for a smooth projective symplectic or Poisson surface,
by the classification of such surfaces \cite{bartocci-macri}.
Consequently, the exterior product 
(\ref{eq-exterior-product-from-cartesian-square-of-moduli-space})
is an isomorphism. 
%Choose a basis $\{x_1, \dots, x_n\}$ of $K^*_{top}(S)$,
%which is a union of bases of
%the summands $K^0_{top}(S)$ and $K^1_{top}(S)$. 
%We get the K\"{u}nneth decomposition
%(\ref{eq-e-i}), of the class 
%$[\E] \in K^0_{top}(S\times \M_H(v))$ of the universal sheaf, 
%where each class $e_i$ is either in $K^0_{top}(\M_H(v))$ or in 
%$K^1_{top}(\M_H(v))$. 
%
%With the above definitions, Proposition
%\ref{prop-integral-generators} has an analogue for moduli spaces of sheaves, 
%over any smooth projective symplectic or Poisson surface $S$.

\begin{prop}
\label{prop-integral-generators-non-simply-connected-surface}
Theorem \ref{thm-introduction-integral-generators} holds, under the 
additional 
assumption that a universal sheaf $\E$ exists over $S\times \M_H(v)$.
%\begin{enumerate}
%\item
%The cohomology ring $H^*(\M_H(v),\Integers)$ is generated by 
%the Chern classes $c_j(e_i)$, of the classes $e_i\in K^*_{top}\M_H(v)$,
%which are given in equation
%(\ref). 
%\item
%The cohomology groups $H^i(\M_H(v),\Integers)$ are torsion free for all $i$. 
%\end{enumerate}
\end{prop}

\noindent
{\bf Proof:}
The proof of Proposition \ref{prop-integral-generators} 
extends to the more general setup. 
In formula (\ref{eq-first-decomposition-of-diagonal}) we implicitly used
the K\"{u}nneth decomposition of 
$\E^\vee=\sum_{i=1}^nx_i^\vee\otimes e_i^\vee$. 
When $x_i$ and $e_i$ 
are classes in $K^1_{top}$, the exterior product 
$f_1^!x_i^\vee\cup f_2^!e_i^\vee$ need not be equal to 
$(f_1^!x_i\cup f_2^!e_i)^\vee$, if we interpret 
$^\vee:K^1_{top}\rightarrow K^1_{top}$ as the duality operator on the 
suspension. Instead, we avoid relating the K\"{u}nneth decompositions 
of $\E$ and $\E^\vee$, write
$[\E^\vee]=\sum_{i=1}^n x_i\otimes e_i'$, replace the Mukai pairing by
$(x,y):=-\chi(x\cup y)$, and replace $e_i^\vee$ by $e_i'$ in 
(\ref{eq-translation-of-thm-diagonal}). 

Somewhat delicate is the analogue of the statement, 
that the classes $\beta_j$, 
in the decomposition (\ref{eq-decomposition-of-the-diagonal-of-moduli}),
are polynomials, with integral coefficients, in the Chern classes of
the $e_i$, for $1\leq i \leq n$. The analogous statement follows from Lemma
\ref{lemma-integral-coefficients}.
\EndProof

\begin{new-lemma}
\label{lemma-integral-coefficients}
Let $x, y\in K^1_{top}X$ and $d$ an integer $\geq 1$. 
The class $c_d(x\cup y)$ in
$H^{2d}(X,\Integers)$ can be written as a polynomial, with 
integral coefficients, in the even-dimensional classes
$c_{i+\frac{1}{2}}(x)\cup c_{k-i-\frac{1}{2}}(y)$, 
for $0\leq i \leq k-1$ and $1\leq k\leq d$.
\end{new-lemma}

The proof of Lemma \ref{lemma-integral-coefficients}
will depend on Lemmas 
\ref{lemma-Chern-character-of-odd-classes} and 
\ref{lemma-reduction-to-torsion-free-case}.
Let $\P_n$, $n\geq 1$, be the set of descending partitions 
$\lambda_1\geq \lambda_2 \geq \cdots \geq \cdots$, $\lambda_i\geq 0$, 
$\sum_{i=1}^\infty\lambda_i=n$. The length $\ell(\lambda)$ is
$\max\{i \ : \ \lambda_i\neq 0\}$. 
The multiplicity $m_i:=m_i(\lambda)$ of $i$ in $\lambda$ is 
the number of $j$ with $\lambda_j=i$.

\begin{new-lemma}
\label{lemma-Chern-character-of-odd-classes}
Let $X$ be a topological space and $x$ a class in $K^*_{top}(X)$. 
Write $x=y+z$, where $y\in K^0_{top}(X)$ and 
$z\in K^1_{top}(X)$. Let $ch_{i}(x)$ be the degree
$2i$ summand  of $ch(x)$ in $H^{2i}(X,\RationalNumbers)$.
Let $k\geq 1$ be an integer.
\begin{enumerate} 
\item
\label{lemma-item-Chern-character-for-even-classes}
${\displaystyle 
(-1)^k(k-1)!ch_k(x)  =
\sum_{i_1+2i_2+ \cdot + ki_k=k}(-1)^{i_1+\cdots +i_k}
\frac{(i_1 + \cdots + i_k-1)!}{i_1!\cdots i_k!}
c_1(y)^{i_1}\cdots c_k(y)^{i_k}.
}$
\item
\label{lemma-item-Chern-character-for-odd-classes}
${\displaystyle ch_{k-\frac{1}{2}}(x)=
\frac{(-1)^{k-1}}{(k-1)!}c_{k-\frac{1}{2}}(z).
}$
\item
\label{lemma-item-Chern-class-in-terms-of-character}
${\displaystyle
c_k(y) = \sum_{\lambda\in \P_k}(-1)^{k-\ell(\lambda)}
\prod_{i\geq 1}\frac{[(i-1)!ch_i(y)]^{m_i}}{m_i!}.
}$
\end{enumerate}
\end{new-lemma}

\noindent 
{\bf Proof:}
\ref{lemma-item-Chern-character-for-even-classes})
The equality $ch_k(x)=ch_k(y)$ follows from the linearity of the 
Chern character. $ch_k(y)$ is a polynomial in 
$c_1(y)$, \dots, $c_{i-1}(y)$. 
It suffices to calculate its coefficients when 
the class $y$ is represented by a vector bundle. For such $y$,
the equality is Girard's formula (\cite{milnor} Ch 16 Problem 16-A).
%The difference $i!ch_i(x)-(-1)^{i-1}ic_i(y)$ is a polynomial in
%$c_1(y)$, \dots, $c_{i-1}(y)$, by the linearity of the Chern character.
%It suffices to prove the integrality of its coefficients when 
%the class $y$ is represented by a vector bundle.
%Using the Splitting Principle,
%$c_i(y)$ is the $i$-th elementary symmetric function in the Chern roots, 
%while $(i!)ch_i(y)$ is the $i$-th power sum of the Chern roots. 
%The statement follows from the determinantal expression 
%of the power sums as polynomials in the elementary symmetric functions
%(\cite{macdonald}, Exercise 8 page 28).

\ref{lemma-item-Chern-character-for-odd-classes})
Let $\tilde{z}\in \widetilde{K}^0_{top}(SX^+)$ be the class corresponding 
to $z$. 
The Chern classes $c_i(\tilde{z})$ belong to the reduced cohomology
$\widetilde{H}^*(SX^+)$, with integral coefficients. 
Let $S^1$ be the $1$-sphere and $u\in H^1(S^1,\Integers)$ a generator.
Choose a base point $s_0\in S^1$. 
Recall the short exact sequence:
\[
0\rightarrow
\widetilde{H}^*(SX^+) \rightarrow
\widetilde{H}^*(S^1\times X) \rightarrow
%\widetilde{H}^*(S^1\times\{\xi_0\})\ \oplus \
\widetilde{H}^*(\{s_0\}\times X) \rightarrow 0
\]
(\cite{atiyah-book}, Corollary 2.4.7). 
In particular, the image of $\widetilde{H}^*(SX^+)$ 
is the principal ideal generated in 
$\widetilde{H}^*(S^1\times X)$ by $u\otimes 1$. 
The cup product $\cup : \widetilde{H}^*(SX^+)\otimes \widetilde{H}^*(SX^+)
\rightarrow \widetilde{H}^*(SX^+)$ vanishes, since $u\cup u=0$ in 
$H^*(S^1)$. Consequently, $c_i(\tilde{z})\cup c_j(\tilde{z})=0$, for any 
two positive integers $i$, $j$. 
Part \ref{lemma-item-Chern-character-for-even-classes} 
of the Lemma implies the equality
\begin{equation}
\label{eq-ch-i-tilde-z-in-terms-of-c-i-tilde-z}
ch_i(\tilde{z}) \ \ = \ \ \frac{(-1)^{i-1}}{(i-1)!}c_i(\tilde{z}).
\end{equation}

Let $\pi_i$ be the projection from $S^1\times X$ to the $i$-th
factor, $i=1,2$. Let $v\in K^1_{top}(S^1)$ be the generator with
$c_{\frac{1}{2}}(v)=u$. Then $\tilde{z}=\pi_1^!(v)\cup\pi_2^!(z)$, 
$ch(\tilde{z})=\pi_1^*ch(v)\cup \pi_2^*ch(z)$, and 
$ch(v)=u$. We get the equality
$
ch_i(\tilde{z})  =  \pi_1^*(u)\cup \pi_2^*ch_{i-\frac{1}{2}}(z).
$
Part \ref{lemma-item-Chern-character-for-odd-classes} of the
Lemma follows from the above equality and equation
(\ref{eq-ch-i-tilde-z-in-terms-of-c-i-tilde-z}).

\ref{lemma-item-Chern-class-in-terms-of-character}) See
\cite{macdonald} Ch. I section 2 equation (2.14').
\EndProof

\begin{new-lemma}
\label{lemma-reduction-to-torsion-free-case}
(\cite{atiyah-kunneth} Lemma 2). 
Let $X$ be a finite CW-complex. Then $X$ can be embedded as a subcomplex 
of a finite CW-complex $A$, so that both $H^*(A,\Integers)$ and
$K^*(A)$ are free abelian groups and 
$K^*(A)\rightarrow K^*(X)$ is surjective.
\end{new-lemma}

\noindent
{\bf Proof of Lemma \ref{lemma-integral-coefficients}:}
If $f:X\rightarrow A$ is continuous,
$x=f^!(x')$, $y=f^!(y')$ and Lemma
\ref{lemma-integral-coefficients} holds for the class
$c_d(x'\cup y')$, then the lemma
holds also for the class $c_d(x\cup y)$. 
We may thus assume that  $H^*(X,\Integers)$ is 
free, by Lemma \ref{lemma-reduction-to-torsion-free-case}.

%The proof is by induction on $d$. The case $d=1$ is clear. 
%Assume that the lemma holds for $k<d$. 
Part \ref{lemma-item-Chern-character-for-odd-classes} of 
Lemma \ref{lemma-Chern-character-of-odd-classes} and the 
multiplicative property of the Chern character yield:
\begin{equation}
\label{eq-ch-d-as-a-sum}
(d-1)!ch_d(x\cup y) \ \ \ = \ \ \ 
(-1)^{d-1}\sum_{i=0}^{d-1}
\Choose{d-1}{i}c_{i+\frac{1}{2}}(x)\cup 
c_{d-i-\frac{1}{2}}(y).
\end{equation}
Hence $(d-1)! ch_d(x\cup y)$ belongs to the subring of
$H^*(X,\Integers)$ generated by 
$c_{i-\frac{1}{2}}(x)\cup c_{j-\frac{1}{2}}(y)$, $1\leq i,j \leq d$. 
The powers $[c_{i-\frac{1}{2}}(x)\cup c_{j-\frac{1}{2}}(y)]^k$ vanish
for $k\geq 2$. We get the identity
\[
[(d-1)!ch_d(x\cup y)]^m = 
m!(-1)^{m(d-1)}\sum_{0\leq k_1 < k_2 < \cdots < k_m \leq d-1}
\prod_{i=1}^m\left[
\Choose{d-1}{k_i}c_{k_i+\frac{1}{2}}(x)\cup c_{d-k_i-\frac{1}{2}}(y)
\right]
\]
Lemma \ref{lemma-integral-coefficients} follows from the above equality and 
part \ref{lemma-item-Chern-class-in-terms-of-character} of 
Lemma \ref{lemma-Chern-character-of-odd-classes}.
\EndProof

\section{Generators in the absence of a universal sheaf}
\label{sec-twisted-sheaves}

We complete the proof of Theorem \ref{thm-introduction-integral-generators}
in this section, dropping the assumption that a universal sheaf exists. 
We construct a universal class $e$ 
in $K^0_{top}(S\times \M)$, regardless of the existence of an 
algebraic universal sheaf (Definition \ref{def-e-v}). 
The details are worked out in section \ref{sec-universal-class}.
In section 
\ref{sec-integral-generators-with-out-univ-sheaf} 
we prove the following Proposition.
Let $S$, $v$, $H$ satisfy the assumptions of Theorem \ref{thm-diagonal},
with one exception: A universal sheaf need not exist. 

\begin{prop}
\label{prop-class-of-diagonal-without-universal-sheaf}
\begin{enumerate}
\item
\label{prop-item-class-of-diagonal-without-universal-sheaf}
Using the notation of Theorem \ref{thm-diagonal}, the class
%Theorem \ref{thm-diagonal} holds for $\M_H(v)$, with the class
%(\ref{eq-class-of-diagonal}) replaced by
\begin{equation}
\label{eq-class-of-diagonal-in-terms-of-C-infinity-universal-sheaf}
c_m\left[- \ 
\pi_{13_!}\left(
\pi_{12}^!(e)^\vee\cup\pi_{23}^!(e)
\right)
\right]
\end{equation}
is Poincare-dual to the class of the diagonal, 
where $e$ is the class in Definition \ref{def-e-v}.
\item
\label{prop-item-derivative-class-vanishes}
When $S$ is a $K3$ or abelian surface, the following vanishing holds
\begin{equation}
\label{eq-derivative-class-in-terms-of-C-infinity-universal-sheaf}
c_{m-1}\left[- \ 
\pi_{13_!}\left(
\pi_{12}^!(e)^\vee\cup\pi_{23}^!(e)
\right)
\right] \ \ = \ \ 0.
\end{equation}
\end{enumerate}
\end{prop}

An immediate consequence is: 

\begin{cor}
\label{cor-generators-in-absence-of-universal-sheaf}
Theorem \ref{thm-introduction-integral-generators}
holds also in the absence of a universal sheaf, 
once we replace in (\ref{introduction-eq-e-i}) the class $\E$, 
of the universal sheaf, by the class 
$e$ in Definition \ref{def-e-v}.
\end{cor}

Following is a summary of the construction of the universal class $e$. 
Let $X$ be a topological space, a complex analytic space, or a scheme 
over $\ComplexNumbers$ endowed with the \'{e}tale topology.
Denote by $\StructureSheaf{X}^*$ the sheaf of  invertible complex valued
continuous, holomorphic, or 
algebraic functions. 
Let $\theta$ be a \v{C}ech $2$-cocycle with coefficients in 
$\StructureSheaf{X}^*$. 
%Let $X$ be a topological space, $\A^*$ (???) the sheaf  of complex valued 
%invertible continuous functions on $X$, 
%and $\theta$ a \v{C}ech $2$-cocycle with coefficient in $\A^*$. 
There is a notion of $\theta$-twisted vector bundles over $X$
(Definition \ref{def-twisted-sheaves}). 
One can define the $K$-group $K^0(X)_{\theta}$ of
$\theta$-twisted vector bundles. The topological version
$K^0_{top}(X)_{\theta}$ is defined in 
\cite{atiyah-segal,donovan-karoubi}. For the analytic 
$K^0_{hol}(X)_{\theta}$ 
or algebraic $K^0_{alg}(X)_{\theta}$ see 
\cite{calduraru-thesis,donagi-pantev}.
$K^0(X)_{\theta}$ is a $K^0(X)$-module. 
$K^0(X)_{\theta}$ 
depends only on the \v{C}ech cohomology class
of $\theta$, canonically up to tensorization by the class
in $K^0(X)$ of a line-bundle. 
Note, that in the topological category 
$H^2(X,\StructureSheaf{X}^*)$ is isomorphic to $H^3(X,\Integers)$, via the 
connecting homomorphism of the exponential sequence. 
%We denote by 
%$[\theta]$ the class corresponding to $\theta$ in $H^3(X,\Integers)$.

Let $S$ be a K3, abelian, or a smooth projective Poisson surface,
and $\M_H(v)$ a moduli space
as in Theorem \ref{thm-introduction-integral-generators}. 
Over $S\times \M_H(v)$ there 
always is a twisted universal sheaf $\E_v$.
The twisting cocycle is the pullback  $f_2^*\theta$ 
of some \v{C}ech $2$-cocycle $\theta$, 
with coefficient in the sheaf $\StructureSheaf{\M_H(v)}^*$,
in the \'{e}tale or classical topology  of $\M_H(v)$ 
(see the Appendix in \cite{mukai-hodge} or section 
\ref{sec-universal-class} below). 
The twisted universal sheaf $\E_v$ 
%can be resolved by $f_2^*\theta$-twisted vector bundles, so it 
determines a class in 
$K^0_{hol}(\M_H(v)\times S)_{f_2^*\theta}$. 
Let $\bar{\theta}$ be the image of $\theta$, as a
\v{C}ech $2$-cocycle with coefficient in the sheaf of 
invertible continuous functions. 
We prove that $\bar{\theta}$ is a coboundary 
(Lemma \ref{lemma-cohomological-triviality-of-theta}).
It follows that 
$K_{top}(\M_H(v)\times S)_{f_2^*\bar{\theta}}$ is isomorphic to 
the untwisted group $K_{top}(\M_H(v)\times S)$,
canonically up to tensorization of the latter with a topological line-bundle
on $\M_H(v)$. 

\begin{defi}
\label{def-e-v}
The universal class $e$ is the image of $\E_v$ under the composition
\[
K_{hol}(\M_H(v)\times S)_{f_2^*\theta}\rightarrow
K_{top}(\M_H(v)\times S)_{f_2^*\bar{\theta}}\rightarrow
K_{top}(\M_H(v)\times S).
\]
\end{defi}

The details are worked out in section \ref{sec-universal-class}.

%***************************************************************************
% A universal class in $K^0_{top}(S\times \M)$
%***************************************************************************
\subsection{A universal class in $K^0_{top}(S\times \M)$}
\label{sec-universal-class}

%We prove in this section, that a universal class 
%(\ref{eq-differentiable-universal-sheaf}) 
%in $K^0_{top}(S\times \M)$ exists, regardless of the existence of an 
%algebraic universal sheaf.
Let $S$ be a K3, abelian, or a smooth projective Poisson surface. 
The moduli space $\M_H(v)$, with a $v$-generic polarization $H$,
always admits a twisted universal sheaf (see the Appendix in 
\cite{mukai-hodge}).
There exists a covering $\U:=\{U_\alpha\}_{\alpha\in I}$ of $\M_H(v)$, 
in the \'{e}tale or classical topology, 
universal sheaves $\E_\alpha$ over $S\times U_\alpha$, isomorphisms
\[
%\begin{equation}
%\label{eq-g-i-j}
g_{\alpha\beta} \ : \ (\E_\beta\restricted{)}{[S\times U_{\alpha\beta}]}
\ \ \LongIsomRightArrow \ \ 
(\E_\alpha\restricted{)}{[S\times U_{\alpha\beta}]},
%\end{equation}
\]
such that the cocycle
\[
(\delta g)_{\alpha\beta\gamma} \ \ := \ \ 
g_{\alpha\beta}g_{\beta\gamma}g_{\gamma\alpha}
\]
comes from a $2$-cocycle $\theta$ in $Z^2(\U,\StructureSheaf{\M_H(v)}^*)$. 

Given a class in $K^0_{alg}(S)$, represented by 
a complex $F$, we get the line-bundle
$L_\alpha:=\det f_{2_!}(\E_\alpha\otimes f_1^*F)$ on $U_\alpha$. The
isomorphisms $g_{\alpha\beta}$ induce 
isomorphisms 
\begin{equation}
\label{eq-g-F-alpha-beta}
g^F_{\alpha\beta}\ : \ (L_\beta\restricted{)}{U_{\alpha\beta}}
\ \ \LongIsomRightArrow \ \ (L_\alpha\restricted{)}{U_{\alpha\beta}},
\end{equation}
satisfying $(\delta g^F)_{\alpha\beta\gamma} = 
\theta^{\chi(v\cup F)}_{\alpha\beta\gamma}$. 
Consequently, $\theta^{\chi(v\cup F)}$ is a coboundary. 
We conclude, that $\theta^n$ is a coboundary, where the natural number $n$ 
is given by 
\begin{equation}
\label{eq-n}
n \ \ := \ \ \mbox{g.c.d}\{
\chi(v\cup w) \ \ : \ \  w\in K^0_{alg}(S)
\}. 
\end{equation}
If $n=1$, then $\theta$ is a coboundary $\theta=\delta\eta$,
for some one-cochain $\eta\in C^1(\U,\StructureSheaf{\M_H(v)}^*)$. 
Then the new transition functions 
\begin{equation}
\label{eq-varphi}
\varphi_{\alpha\beta} \ := g_{\alpha\beta} \cdot \eta_{\alpha\beta}^{-1}
\end{equation}
glue the local universal sheaves $\E_\alpha$ to a global universal sheaf. 
In general, we get only a $f_2^*\theta$-twisted sheaf
$\E_v$, in the following sense.

\begin{defi}
\label{def-twisted-sheaves}
{\rm
Let $X$ be a scheme or a complex analytic space,
$\U:=\{U_\alpha\}_{\alpha\in I}$ a covering, open in the complex 
or \'{e}tale topology,
and $\theta\in Z^2(\U,\StructureSheaf{X}^*)$ a \v{C}ech $2$-cocycle.
A {\em $\theta$-twisted sheaf} consists of sheaves $\E_\alpha$ of
$\StructureSheaf{U_\alpha}$-modules over $U_\alpha$,
for all $\alpha\in I$, and isomorphisms 
$g_{\alpha\beta}:(\E_\beta\restricted{)}{U_{\alpha\beta}}
\rightarrow (\E_\alpha\restricted{)}{U_{\alpha\beta}}$ 
%as in (\ref{eq-g-i-j}) 
satisfying the conditions:

(1) $g_{\alpha\alpha}=id$,

(2) $g_{\alpha\beta}=g_{\beta\alpha}^{-1}$, 

(3) $g_{\alpha\beta}g_{\beta\gamma}g_{\gamma\alpha}=
\theta_{\alpha\beta\gamma}\cdot id.$ 

\noindent
The $\theta$-twisted sheaf is {\em coherent}, if the $\E_i$ are.
}
\end{defi}

Locally free $\theta$-twisted sheaves of finite rank form an abelian category,
with the obvious notion of homomorphisms, and we let 
$K^0_{hol}(X)_\theta$ be its $K$-group. 
Observe, that the determinant $\det(\E)$, of a $\theta$-twisted locally free 
sheaf $\E$ of rank $r$, is a $\theta^r$-twisted line-bundle.
Thus, $\theta^r$ is a coboundary. Consequently, 
the order of the class $[\theta]$,
of $\theta$ in $H^2(X,\StructureSheaf{X}^*)$,
divides the rank of every $\theta$-twisted locally free sheaf $\E$. 
$K^0_{hol}(X)_\theta$ is thus trivial, if $[\theta]$ has infinite order. 
For a more general definition of twisted $K$-groups and 
derived categories, see \cite{donagi-pantev} and references therein.

Let $\theta$ be a \v{C}ech $2$-cocycle of continuous invertible
complex valued functions. 
The twisted topological $K$-groups are defined 
for a class $[\theta]$ of arbitrary order \cite{atiyah-segal}.
When the order of of  $[\theta]$ 
is finite, the obvious analogue of the above definition
yields a group $K^0_{top}(X)_\theta$, which is naturally isomorphic
to the one defined in \cite{atiyah-segal}.

We define the class $[\E_v]$
of the twisted universal sheaf 
$\E_v$ in $K_{hol}^0(S\times \M_H(v))_{f_2^*\theta}$
using the following twisted locally free resolution of $\E_v$. 
Choose a sufficiently ample line-bundle
$H'$ on $S$ and set $W_\alpha:=f_{2_*}(\E_\alpha\otimes f_1^*H')$ and 
$\E_{\alpha,0}:= f_1^*(H')^{-1}\otimes f_2^*W_\alpha$. 
There is a natural evaluation homomorphism 
$ev_\alpha:\E_{\alpha,0}\rightarrow \E_\alpha$.
We choose $H'$ sufficiently ample, so that 
$R^{i}f_{2_*}(\E_\alpha\otimes f_1^*H')$ vanishes, for $i>0$, 
each $\E_{\alpha,0}$ is locally free, and each $ev_\alpha$ is surjective. 
The gluing transformations $g_{\alpha\beta}$, of the $\E_\alpha$, induce 
gluing transformations $\psi_{\alpha\beta}$ 
of the vector bundles $W_\alpha$, 
whose coboundary is diagonal $\theta_{\alpha\beta\gamma}\cdot id$. 
We denote the corresponding $\theta$-twisted locally free sheaf by 
\begin{equation}
\label{eq-twisted-vector-bundle-W}
W.
\end{equation}
Similarly, the pullback $f_2^*(\psi_{\alpha\beta})$ 
is a $1$-cochain of gluing transformations for the $\E_{\alpha,0}$,
defining a $f_2^*\theta$-twisted vector bundle $\E_0$.
Repeating the process once more, starting with $\E_0$, we get 
a surjective homomorphism $\E_1\rightarrow \E_0$,
with kernel $\E_2$, where $\E_i$, $0\leq i \leq 2$, are 
locally free $f_2^*\theta$-twisted sheaves. 
The class $\sum(-1)^i[\E_i]$ is independent of the choices 
made and is denoted by $[\E_v]$.

When the cocycle $\theta$ is a coboundary,
$\theta=\delta(\eta)$, the procedure described in equation 
(\ref{eq-varphi}) induces a well defined isomorphism
$K_{hol}^0(X)_\theta\rightarrow K_{hol}^0(X)$.

Let $\A^*$ be the sheaf of continuous, complex valued, and 
invertible functions on $\M_H(v)$. 

\begin{new-lemma}
\label{lemma-cohomological-triviality-of-theta}
The homomorphism $Z^2(\U,\StructureSheaf{\M_H(v)}^*)\rightarrow Z^2(\U,\A^*)$
maps the cocycle $\theta$ to a coboundary $\bar{\theta}$. 
%\begin{enumerate}
%\item
%\label{lemma-item-theta-is-topologically-trivial}
%The homomorphism $Z^2(\U,\StructureSheaf{\M_H(v)}^*)\rightarrow Z^2(\U,\A^*)$
%maps the cocycle $\theta$ to a coboundary $\bar{\theta}$. 
%\item
%\label{lemma-item-phi-is-injective}
%The homomorphism 
%$\phi^*:H^*(\M_H(v),\Integers)\rightarrow H^*(\PP,\Integers)$ 
%is injective.
%\end{enumerate} 
\end{new-lemma}

The general idea for the proof of Lemma 
\ref{lemma-cohomological-triviality-of-theta}
is clear; the number $n$ in equation
(\ref{eq-n}) becomes $1$, once we replace $K^0_{alg}S$ by 
$K^0_{top}S$. Assume the existence of the Gysin map 
$f_{2_!}:K^0_{top}(S\times \M_H(v))_{f_2^*\bar{\theta}}
\rightarrow K^0_{top}(\M_H(v))_{\bar{\theta}}$
in twisted 
topological $K$-theory, 
and choose $x\in K^0_{top}S$ satisfying $\chi(x\cup v)=1$.
Then $f_{2_!}(f_1^!(x)\cup[\E_v])$ is a class 
in $K_{top}^0(\M_H(v))_{\bar{\theta}}$ of rank $1$. 
But the order of the cohomology class of $\bar{\theta}$
divides the rank of any class in $K_{top}^0(\M_H(v))_{\bar{\theta}}$.
Hence $\bar{\theta}$ is a coboundary.
We avoid using Gysin maps in a twisted version of topological
$K$-Theory, as we are unfamiliar with such a construction in the literature. 
The following elementary lemma provides an alternative proof of the 
triviality of $\bar{\theta}$. 
The lemma summarizes  well known facts about 
twisted sheaves, and will be needed also in the next section.

\begin{new-lemma}
\label{lemma-brouer-severi-presentation}
Let $E:=(E_\alpha,g_{\alpha\beta})$ be a $\theta$-twisted sheaf and 
$F:=(F_\alpha,\psi_{\alpha\beta})$ a $\theta$-twisted locally 
free sheaf over an analytic space $X$. 
\begin{enumerate}
\item
\label{lemma-item-projective-bundle}
The projective bundles $\PP{F_\alpha}$ glue to a global projective bundle
$\phi:\PP{F}\rightarrow X$. 
\item
\label{lemma-item-topological-obstruction}
The image of $\theta$ in $H^3(X,\Integers)$, under the 
connecting homomorphism of the exponential sequence, is the topological 
obstruction for $\PP{F}$ to lift to a complex topological 
vector bundle.
\item
\label{lemma-item-L}
There exists a line bundle $L$ over $\PP{F}\times_X\PP{F}$,
which restricts as $\StructureSheaf{}(1,-1)$ to each fiber
$\PP{F_x}\times \PP{F_x}$, over a point $x\in X$.
\item
\label{lemma-item-characterization-of-tilde-E}
There exists a sheaf $\widetilde{E}$ over $\PP{F}$, which restricts 
to $\PP{F_\alpha}$ as 
$\phi^*E_\alpha\otimes \StructureSheaf{\PP{F_\alpha}}(-1)$,
satisfying
\[
p_1^*(\widetilde{E})\otimes L \ \ = \ \ p_2^*\widetilde{E},
\]
where $p_i:\PP{F}\times_X \PP{F}\rightarrow \PP{F}$, $i=1,2$, are the 
projections.
\item
\label{lemma-item-K-class-of-tilde-E-versus-K-class-of-E}
The cocycle $\phi^*\theta$ is a coboundary,
and the class of $\widetilde{E}$ in
$K_{hol}^0(\PP{F})$ is the image of $E$ via the isomorphism
$K_{hol}^0(\PP{F})_{\phi^*\theta}\rightarrow K_{hol}^0(\PP{F})$,
determined by a suitable choice of a $1$-cochain $\tilde{\psi}$
satisfying $\delta(\tilde{\psi})=\phi^*\theta$.
\end{enumerate}
\end{new-lemma}

\noindent
{\bf Proof:}
Part \ref{lemma-item-projective-bundle} is clear and
part \ref{lemma-item-topological-obstruction} is 
standard \cite{calduraru-thesis}. 

Proof of parts 
\ref{lemma-item-L},
\ref{lemma-item-characterization-of-tilde-E}, and 
\ref{lemma-item-K-class-of-tilde-E-versus-K-class-of-E})
Let $\tau_\alpha:=\StructureSheaf{\PP{F_\alpha}}(-1)$ 
be the tautological line sub-bundle over
$\PP{F_\alpha}$. The pullback 
$\phi^*\psi_{\alpha\beta}$ to $\PP{F}$ restricts to an isomorphism 
\[
\tilde{\psi}_{\alpha\beta} \ : \ 
(\tau_\beta\restricted{)}{\PP{F_\alpha}\cap \PP{F_\beta}} \ \ \ 
\longrightarrow \ \ \ 
(\tau_\alpha\restricted{)}{\PP{F_\alpha}\cap \PP{F_\beta}},
\]
defining a $1$-cochain $\tilde{\psi}$ satisfying 
$\delta(\tilde{\psi})=\phi^*\theta$. 
The local sheaves $\phi^*E_\alpha\otimes \tau_\alpha^{-1}$,
over $\PP{F_\alpha}$, are glued to a global sheaf $\widetilde{E}$ 
over $\PP{F}$ via the transformations 
$\phi^*g_{\alpha\beta}\otimes\tilde{\psi}_{\alpha\beta}^{-1}$.
The cocycle 
$p_1^*\tilde{\psi}_{\alpha\beta}\otimes p_2^*\tilde{\psi}_{\alpha\beta}^{-1}$
over $\PP{F}\times_X\PP{F}$ glues 
$p_1^*\tau_\alpha\otimes p_2^*\tau_\alpha^{-1}$ to the global line-bundle $L$.
\EndProof

Let 
\begin{equation}
\label{eq-PP}
\phi:\PP\rightarrow \M_H(v)
\end{equation} 
be a projective bundle, 
corresponding to a $\theta$-twisted vector bundle, such as $W$ 
given in (\ref{eq-twisted-vector-bundle-W}). Let 
$\widetilde{\E}$ be the sheaf over $S\times \PP$, corresponding to the 
twisted universal sheaf $\E_v$, 
via Lemma \ref{lemma-brouer-severi-presentation}. Denote by 
$\tilde{f}_i$  the projection from $S\times \PP$ on the $i$-th factor,
$i=1,2$. 

\medskip
\noindent
{\bf Proof of Lemma \ref{lemma-cohomological-triviality-of-theta}:}
%\ref{lemma-item-theta-is-topologically-trivial}) 
It suffices to prove that $\PP$ is the projectivization of a topological 
vector bundle, by Lemma \ref{lemma-brouer-severi-presentation} part 
\ref{lemma-item-topological-obstruction}.
This, in turn, is equivalent to the surjectivity of
the restriction homomorphism 
$H^2(\PP,\Integers)\rightarrow H^2(\PP_m,\Integers)$, where 
$\PP_m$ is a fiber of $\PP$ over a point $m$ in $\M_H(v)$,
by a well known criterion\footnote{
Details are included in the proof 
of Lemma 15 part 1 in the eprint version math.AG/0406016 v1
of this paper.}. 

A class $x\in K_{top}^0(S)$, satisfying $\chi(x\cup v)=1$,
exists by Remark \ref{rem-perfect-pairing}, since $v$ is primitive.
Set $y:=\tilde{f}_{2_!}[\tilde{f}_1^!(x)\cup \widetilde{\E}]$.
Then $y$ is a class of rank $1$ in $K_{top}^0(\PP)$ satisfying
the equality
\[
p_1^!(y)\cup \ell \ \ = \ \ p_2^!(y),
\]
where $\ell\in K^0_{top}(\PP\times\PP)$ is the class of 
the line bundle $L$ in Lemma
\ref{lemma-brouer-severi-presentation}, and
$p_i$, $i=1,2$, are the projections from $\PP\times \PP$.
The above equality is verified via the following sequence of 
simpler equalities. Let $\pi_{ij}$ be the projections from
$\PP\times \PP\times S$. 
\[
%p_1^!(y)\cup \ell 
%= 
p_1^!\left(
\tilde{f}_{2_!}[\tilde{f}_1^!(x)\cup \widetilde{\E}]
\right)\cup\ell 
=
\pi_{12_!}\left(
\pi_3^!(x)\cup\pi_{13}^!(\widetilde{\E})\cup\pi_{12}^!\ell
\right)
=
\pi_{12_!}\left(
\pi_3^!(x)\cup \pi_{23}^!(\widetilde{\E})
\right)
=
p_2^!(y).
\]
%The definition of $y$ implies the first equality, 
The last equality follows from property
(\ref{eq-cohomology-and-base-change}) for Gysin maps and for the 
first use also the projection formula
(\ref{eq-projection-formula}).
The second equality follows from part
\ref{lemma-item-characterization-of-tilde-E} of Lemma
\ref{lemma-brouer-severi-presentation}.
We conclude, that $c_1(y)$ restricts to each fiber
$\PP_m$ as a generator of $H^2(\PP_m,\Integers)$.
%This completes the proof of Lemma 
%\ref{lemma-cohomological-triviality-of-theta}.
\EndProof

\medskip
As a consequence of Lemma \ref{lemma-cohomological-triviality-of-theta}, 
we can choose a $1$-cochain $\eta\in C^1(\U,\A^*)$, satisfying
$\theta=\delta\eta$. Modifying the transition functions $g_{\alpha\beta}$,
of each $f_2^*\bar{\theta}$-twisted vector bundle, as in equation 
(\ref{eq-varphi}),
we get an
isomorphism $K_{top}^0(S\times \M_H(v))_{f_2^*\bar{\theta}}\rightarrow 
K_{top}^0(S\times \M_H(v))$. Use this isomorphism in 
Definition \ref{def-e-v} to map the class of 
$\E_v$ in $K_{hol}^0(S\times \M_H(v))_{f_2^*\theta}$ to
a  class 
\begin{equation}
\label{eq-differentiable-universal-sheaf}
e \ \ \in \ \ K^0_{top}[S\times \M_H(v)]. 
\end{equation}

%*****************************************************************
% Integral generators
%*****************************************************************
\subsection{Decomposition of the diagonal via a universal class}
\label{sec-integral-generators-with-out-univ-sheaf}

\noindent
{\bf Proof of Proposition 
\ref{prop-class-of-diagonal-without-universal-sheaf}:} 
\ref{prop-item-class-of-diagonal-without-universal-sheaf})
We will prove the equality of the class
(\ref{eq-class-of-diagonal-in-terms-of-C-infinity-universal-sheaf}) 
and the class of the diagonal, by pulling back both via the
homomorphism 
\begin{equation}
\label{eq-phi-times-phi}
(\phi\times\phi)^*:H^*(\M_H(v)\times\M_H(v),\Integers)
\rightarrow H^*(\PP\times \PP,\Integers),
\end{equation}
and comparing both pullbacks to a third class 
(\ref{eq-pulled-back-diagonal}).
$\PP$ is the projectivization of a global topological
vector bundle, by Lemma 
\ref{lemma-cohomological-triviality-of-theta}. 
The injectivity of 
$\phi^*:H^*(\M_H(v),\Integers) \rightarrow H^*(\PP,\Integers)$,
and hence of (\ref{eq-phi-times-phi}), follows
(\cite{karoubi} Proposition V.3.12). 

An extension of Theorem \ref{thm-diagonal}, carried out  
in \cite{markman-diagonal} section 3, 
implies that the pullback by $\phi\times \phi$,
of the diagonal in $\M_H(v)\times \M_H(v)$, is Poincare dual to the class
\begin{equation}
\label{eq-pulled-back-diagonal}
c_m\left[- \ 
\tilde{\pi}_{13_!}\left(
\tilde{\pi}_{12}^!(\widetilde{\E})^\vee\cup\tilde{\pi}_{23}^!(\widetilde{\E})
\right)
\right],
\end{equation}
where $\tilde{\pi}_{ij}$ is the projection from $\PP\times S\times \PP$
onto the product of the $i$-th and $j$-th factors
Furthermore, the vanishing  
\begin{equation}
\label{eq-derivative-class-via-tilde-E-vanishes}
c_{m-1}\left[- \ 
\tilde{\pi}_{13_!}\left(
\tilde{\pi}_{12}^!(\widetilde{\E})^\vee\cup\tilde{\pi}_{23}^!(\widetilde{\E})
\right)
\right] \ \ = \ \ 0,
\end{equation}
holds when $S$ is a $K3$ or abelian surface 
(equation (\ref{eq-vanishing-ofc-m-1})). 

The class $[\widetilde{\E}]$  of $\widetilde{\E}$ in
$K_{top}^0(S\times \PP)$ is the image of the class
$(id_S\times \phi)^!\E_v$ via the composite homomorphism
$
K_{hol}^0(S\times \PP)_{\phi^*(\theta)}\rightarrow 
K_{top}^0(S\times \PP)_{\phi^*(\bar{\theta})}\rightarrow 
K_{top}^0(S\times \PP), 
$
by Lemma 
\ref{lemma-brouer-severi-presentation} part
\ref{lemma-item-K-class-of-tilde-E-versus-K-class-of-E}.
The same holds for the class $(id_S\times \phi)^!e$, by its definition
(\ref{eq-differentiable-universal-sheaf}).
Hence, there exists a topological line bundle  $F$  over $\PP$,
satisfying
\begin{equation}
\label{eq-relation-between-two-universal-objects}
[\widetilde{\E}] \ \ = \ \ F\cup(id_S\times \phi)^!e.
\end{equation}

We claim, that $\phi\times \phi$ pulls back the class
(\ref{eq-class-of-diagonal-in-terms-of-C-infinity-universal-sheaf})
to the class (\ref{eq-pulled-back-diagonal}). 
The proposition follows from this claim, since the homomorphism 
(\ref{eq-phi-times-phi}) 
is injective. 
The claim follows from 
the equivalence (\ref{eq-relation-between-two-universal-objects}) and 
the invariance of the class 
(\ref{eq-pulled-back-diagonal}), under replacement of $\widetilde{\E}$
by $\widetilde{\E}\otimes F^{-1}$. 
The invariance follows from equation 
(\ref{eq-c-r-plus-2-depends-linearly-on-F}) in Lemma 
\ref{lemma-invariance-of-c-m-under-tensorization},
the vanishing (\ref{eq-derivative-class-via-tilde-E-vanishes})
when the surface $S$ is K3 or abelian,
and the following rank calculation. 

The class 
$-\tilde{\pi}_{13_!}\left(
\tilde{\pi}_{12}^!(\widetilde{\E})^\vee\cup\tilde{\pi}_{23}^!(\widetilde{\E})
\right)$  
in $K^0_{alg}(\PP\times \PP)$ has rank $-\chi(v^\vee\cup v)$. 
%where $E$ is a stable sheaf parametrized by $\M_H(v)$. 
The dimension $m$ of $\M_H(v)$ is either $2-\chi(v^\vee\cup v)$
or $1-\chi(v^\vee\cup v)$, depending on $S$ being symplectic, 
or non-symplectic Poisson. 

\ref{prop-item-derivative-class-vanishes})
Equation  (\ref{eq-derivative-class-in-terms-of-C-infinity-universal-sheaf})
pulls back to equation 
(\ref{eq-derivative-class-via-tilde-E-vanishes}), as seen by 
the equivalence (\ref{eq-relation-between-two-universal-objects}) and 
equation (\ref{eq-c-r-plus-1-is-invariant}) in 
Lemma \ref{lemma-invariance-of-c-m-under-tensorization}.
\EndProof

%****************************************************************
%
%****************************************************************
\section{Higgs bundles}
\label{sec-higgs}

%The Theorem sharpens results of Hausel-Thaddeus 
%and the author, where the cohomology was considered with 
%rational coefficients \cite{ht,markman-diagonal}.

We sketch the proof of Theorem \ref{thm-higgs} in this section. 
Let us first review the geometric set-up of the proof of Theorem 7 
in \cite{markman-diagonal}. Set $\H:=\H_\Sigma(r,d,D)$ and
otherwise keep the notation introduced in the
paragraph preceding Theorem \ref{thm-higgs}.
Let $S$ be the projectivization of the rank $2$ vector bundle 
$K_\Sigma(D)\oplus \StructureSheaf{\Sigma}$, and
$b:S\rightarrow \Sigma$ the bundle map. $S$ is a Poisson surface. 
$\H$ is a Zariski open subset of a compact moduli space
$\M$, of stable sheaves of pure dimension $1$ on $S$. 
The moduli space $\M$ may be singular outside $\H$, as
points in the boundary correspond to sheaves on $S$, which may not satisfy 
Condition \ref{cond-proper-subsheaf}. There exists however a resolution 
$\nu:\widetilde{\M} \rightarrow \M$, which is an isomorphism over
$\H$. Furthermore, $\widetilde{\M}$ is projective and 
the restriction homomorphism from 
$H^*(\widetilde{\M},\Integers)$ to
$H^*(\H,\Integers)$ is surjective. 
Let $\pi_{ij}$ be the projection from
$\H\times S\times \widetilde{\M}$ onto the product of the
$i$-th and $j$-th factors. 
A universal sheaf $\F'$ exists over $S\times \M$ and 
we denote by $\F$ its restriction to $S\times \H$. 
The equality $(b\times id)_!(\F)=[\E]$,
of classes  in $K^0_{alg}(\Sigma\times \H)$,
relates the universal sheaf to the universal bundle. 
The class, in
Borel-Moore homology, of the graph of the embedding 
$\iota:\H\hookrightarrow \widetilde{\M}$ is given by 
$c_m\left[-\pi_{13_!}\left(
\pi_{12}^!(\F)^\vee\cup\pi_{23}^!(id\times\nu)^!(\F')
\right)
\right].$

{\bf Step 1:}
We claim, that the Chern classes of the K\"{u}nneth
factors of $\F$ generate the integral cohomology ring of $\H$. 
This follows from the proof of Proposition 
\ref{prop-integral-generators-non-simply-connected-surface},
modulo minor modifications due to the non-compactness of $\H$.  
One needs to replace homology groups and classes by their 
Borel-Moore analogues, and the Poincare-Duality isomorphism,
by the isomorphism between Borel-Moore homology and singular cohomology
of the smooth varieties.
In equation 
(\ref{eq-translation-of-thm-diagonal}), for example, 
$\delta$ is the singular cohomology class, 
corresponding to the class, in Borel-Moore homology, 
of the graph of the embedding 
$\iota:\H\hookrightarrow \widetilde{\M}$. 
Equation (\ref{eq-decomposition-of-the-diagonal-of-moduli}) then holds,
where $p_i$ are the projections from
$\H\times \widetilde{\M}$, and each
$\alpha_j$ (respectively $\beta_j$)
is a super-symmetric polynomial, with integral coefficients, 
in the Chern classes of the K\"{u}nneth factors of $\F$
(respectively $(id\times \nu)^!\F'$). 
Here we use the fact 
that $p_1$ is smooth and proper, which is the
reason for the construction of the auxiliary space $\widetilde{\M}$. 
Equations 
(\ref{eq-gysin-composed-with-cup-delta-is-identity}) and 
(\ref{eq-x-is-a-linear-combination-of-alpha-j})
are replaced by 
\begin{eqnarray}
\nonumber
%\label{eq-gysin-composed-with-cup-delta-is-pull-back-by-iota}
\iota^*(x) & = & p_{1,*}(\delta\cup p_2^*(x)), \ \ \ \mbox{and}
\\
\label{eq-iota-pulls-back-x-to-a-linear-combination-of-alpha-j}
\iota^*(x) & = & 
\sum_{j\in J} \left(\int_{\widetilde{\M}} 
x\cup\beta_j\right)\alpha_j,
\end{eqnarray}
for $x\in H^*(\widetilde{\M},\Integers)$. 
We conclude that  the Chern classes of the K\"{u}nneth factors of 
$\F$ generate $H^*(\H,\Integers)$, by 
equation 
(\ref{eq-iota-pulls-back-x-to-a-linear-combination-of-alpha-j})
and the surjectivity of
$\iota^*: H^*(\widetilde{\M},\Integers)\rightarrow H^*(\H,\Integers)$.

{\bf Step 2:}
It remains to relate the K\"{u}nneth factors of $\F$ to those of $\E$. 
Let $h$ be the class in $K_{top}^*(S)$ of the
line bundle $\StructureSheaf{S}(-1)$. Then 
$K_{top}^*(S)$ is a free $K_{top}^*(\Sigma)$-module with basis
$\{1,h\}$, by \cite{karoubi} Chapter IV, Theorem 2.16.
Let $\tilde{f}_i$, $i=1,2$, be the projections from 
$S\times \H$ and $f_i$, $i=1,2$, the projections from 
$\Sigma\times \H$.
The following two equalities hold for all $x\in K_{top}^*(\Sigma)$.
\begin{eqnarray}
\label{eq-vanishing-kunneth-factors}
\tilde{f}_{2_!}\left(
\tilde{f}_1^![(h-1)b^!(x)]\cup \F
\right) & = & 0,
\\
\label{eq-same-kunneth-factors}
\tilde{f}_{2_!}\left(
\tilde{f}_1^!b^!(x)\cup \F
\right) & = & f_{2_!}(f_1^!(x)\cup \E).
\end{eqnarray}
Equality (\ref{eq-vanishing-kunneth-factors}) follows from the fact,
that the support of $\F$ is the universal spectral curve, which 
is contained in the open subset $K_\Sigma(D)\times \H$ of $S\times \H$, 
where $K_\Sigma(D)$ denotes also the total space of the line-bundle. 
Equality (\ref{eq-same-kunneth-factors}) follows from
the relation $(b\times id)_!(\F)=\E$ and the projection
formula (\ref{eq-projection-formula}). 
Equalities (\ref{eq-vanishing-kunneth-factors}) and 
(\ref{eq-same-kunneth-factors}) imply, that the K\"{u}nneth factors of 
$\F$ span the same subgroup of
$K_{top}^*(\H)$ as those of $\E$. 
\EndProof

%****************************************************************
% Bibliography:
%****************************************************************

University of Massachusetts, Amherst, MA 01003, USA  

E-mail: markman@math.umass.edu

\end{document}